\documentclass[11pt]{article}
\usepackage{epsfig}
\usepackage{amsfonts}

\newenvironment{proof}{\noindent{\bf Proof.\,}}{\hfill$\Box$}
\newenvironment{proofcite}[1]{\noindent{\bf Proof of #1.\,}}{\hfill$\Box$}

\newtheorem{theorem}{Theorem}[section]
\newtheorem{lemma}[theorem]{Lemma}
\newtheorem{corollary}[theorem]{Corollary}
\newtheorem{proposition}[theorem]{Proposition}
\newtheorem{definition}[theorem]{Definition}
\newtheorem{conjecture}[theorem]{Conjecture}

\def\eps{\varepsilon}
\def\ie{{\em i.e.}}

\newcommand{\subsubsubsection}[1]{\begin{center} {\bf
#1}\end{center}}

\input epsf

\begin{document}

\title{Tiling tripartite graphs with $3$-colorable graphs}
\author{{\bf Ryan R. Martin}\thanks{Corresponding author.
Research supported in part by
NSA grant H98230-05-1-0257. Email: {\tt rymartin@iastate.edu}}
\\ Iowa State University \\ Ames, IA 50010 \and {\bf Yi
Zhao}\thanks{Research supported in part by NSA grant
H98230-05-1-0079. Part of this research was done while working at
University of Illinois
at Chicago. Email: {\tt matyxz@langate.gsu.edu}} \\
Georgia State University \\ Atlanta, GA 30303}
%\date{}

\maketitle

\begin{abstract}
   For a fixed integer $h\geq 1$, let $G$ be a tripartite graph with
   $N$ vertices in each vertex class, $N$ divisible by $6h$, such
   that every vertex is adjacent to at least
   $2N/3+h-1$ vertices
   in each of the other classes. We show that if $N$ is
   sufficiently large, then $G$ can be tiled perfectly by copies
   of $K_{h,h,h}$.  This extends the work in~\cite{MM} and also
   gives a sufficient condition for tiling by any (fixed)
   3-colorable graph.  Furthermore, we show that this
   minimum-degree condition is best possible and provide very tight bounds when $N$ is divisible by $h$ but not by $6h$.
\end{abstract}

\section{Introduction}
Let $H$ be a graph on $h$  vertices, and let $G$ be a graph on $n$
vertices.  Tiling problems in extremal graph theory are
investigations of the condition or conditions under which $G$ must
contain many vertex disjoint copies of $H$ (as subgraphs). An $H$-tiling of $G$ is a subgraph of $G$ which consists of vertex-disjoint copies of $H$.  A {\em perfect $H$-tiling} of $G$ is an $H$-tiling
consisting of $\lfloor n/h\rfloor$ copies of $H$.  For clarity and consistency with other results in this area, we call a perfect $H$-tiling an {\em $H$-factor}.  A very early
tiling result is Dirac's theorem on Hamilton cycles~\cite{Dirac},
which implies that every $n$-vertex graph $G$ with minimum degree
$\delta(G)\ge n/2$ contains a perfect matching (usually called
1-factor, instead of $K_2$-factor).  Later Corr\'adi and
Hajnal~\cite{CoHa} studied the minimum degree of $G$ that guarantees
a $K_3$-factor.  Hajnal and Szemer\'edi \cite{HaSz} settled the
tiling problem for any complete graph $K_h$ by showing that every
$n$-vertex graph $G$ with $\delta(G)\ge (h-1)n/h$ contains a
$K_r$-factor (it is easy to see that this is sharp).  Using the
celebrated Regularity Lemma of Szemer\'edi \cite{Sz}, Alon and
Yuster \cite{AlYu1, AlYu2} obtained results on $H$-tiling for
arbitrary $H$. Their results were later improved by various
researchers \cite{KSSz-AY, Komlos, AliYi, KuhnOsthus}.

In this paper, we consider multipartite tiling, which restricts $G$
to be an $r$-partite graph. When $r=2$, The K\"{o}nig-Hall Theorem
(e.g. see \cite{Bollobas}) answers the 1-factor problem for
bipartite graphs. Wang \cite{Wang98, Wang99} considered
$K_{s,s}$-factors in bipartite graphs for all $s>1$, the second
author \cite{Zhao} gave the best possible minimum degree condition
for this problem.

In a tripartite graph $G=(A,B,C;E)$, the graphs induced by $(A,B)$,
$(A,C)$ and $(B,C)$ are called the \textbf{natural bipartite
subgraphs} of $G$.  Let ${\cal G}_r(N)$ denote the family of
$r$-partite graphs with $N$ vertices in each of its partition sets.
In an $r$-partite graph $G$, $\bar{\delta}(G)$ stands for the
minimum degree from a vertex in one partition set to any other
partition set. Fischer~\cite{Fischer} gives almost perfect
$K_3$-tilings in ${\cal G}_3(N)$ with $\bar{\delta}(G)\geq 2N/3$ and
Johansson~\cite{Johansson} gives a $K_3$-factor with a less
stringent degree condition $\bar{\delta}(G)\geq
2N/3+O\left(\sqrt{N}\right)$. For all $r>2$, Fischer \cite{Fischer}
conjectured the following variant of Hajnal-Szemer\'edi Theorem.
\begin{conjecture}[Fischer \cite{Fischer}]
\label{conj1} If $G\in {\cal G}_r(N)$ satisfies $\bar{\delta}(G)\ge
\frac{r-1}r N$, then $G$ contains a $K_r$-factor.
\end{conjecture}

Recently, Szemer\'edi and the first author \cite{MSz} proved
Conjecture~\ref{conj1} for $r=4$. However, Conjecture~\ref{conj1} is
false when $r=3$: the following construction of Magyar and the first
author, \cite{MM}, provides a counterexample. Let $\Gamma_r\in {\cal
G}_r(r)$ have vertices $\{h_{i,j} : i=1,\ldots,r; j=1,\ldots,r\}$
and the adjacency rule as follows: $h_{i,j}\sim h_{i',j'}$ if $i\neq
i'$ and $j\neq j'$ and either $j$ or $j'$ is in $\{1,\ldots,r-2\}$.
Also, $h_{i,r-1}\sim h_{i',r-1}$ and $h_{i,r}\sim h_{i',r}$ for
$i\neq i'$. No other edges exist. It is easy to check that
$\bar{\delta}(\Gamma_r) = r-1$ and when $r$ is odd, $\Gamma_r$
contains no $K_r$-factor.

Nevertheless, \cite{MM} showed that, if $N$ is an odd multiple of
3, the so-called blow-up graph $\Gamma_3(N)\in {\cal G}_3(N)$
(where each edge of $\Gamma_3$ is replaced with a $K_{N/3,N/3}$
and each non-edge is replaced by an $(N/3)\times (N/3)$ bipartite
graph with no edges) is the unique exception for
Conjecture~\ref{conj1} in the case $r=3$. As a result, this gives
the following Corr\'adi-Hajnal-type theorem.

\begin{theorem}[Magyar-M.~\cite{MM}]
If $G\in {\cal G}_3(N)$ satisfies $\bar{\delta}(G)\ge (2/3)N+1$,
then $G$ contains a $K_3$-factor. \label{thm:MM}
\end{theorem}

In this paper we extend this result to all 3-colorable graphs. Our
main result is the following theorem.

\begin{theorem}
Fix a positive integer $h$.  Let $f(h)$ be the smallest value for which there exists an $N_0$
such that if $G$ is a balanced tripartite graph on $3N$ vertices,
$N\geq N_0$, $h\mid N$, and each vertex is adjacent to at
least $f(h)$ vertices in each of the other classes, then $G$
contains a $K_{h,h,h}$-factor.

If $N=(6q+r)h$ with $0\leq r<6$, then
$$ \begin{array}{lcccll}
   & & f(h) & = & \frac{2N}{3}+h-1, & \mbox{ if $r=0$;} \\
   h\left\lceil\frac{2N}{3h}\right\rceil+h-2
   & \leq & f(h) & \leq &
   h\left\lceil\frac{2N}{3h}\right\rceil+h-1,
   & \mbox{ if $r=1,2,4,5$;} \\
   \frac{2N}{3}+h-1 & \leq & f(h) & \leq & \frac{2N}{3}+2h-1, &
   \mbox{ if $r=3$.}
   \end{array} $$

% and $N$ is divisible by $h$, the following occurs for any
% $G\in {\cal G}_3(N)$:
%    \begin{itemize}
%    \item If each vertex is adjacent to at least $(2/3)N+2h-1$
%    vertices in each of the other classes, then $G$ contains a
%    $K_{h,h,h}$-tiling.
%    \item If $N$ is a multiple of $6h$, and each vertex is
%    adjacent to at least $(2/3)N+h-1$ vertices in each of the
%    other classes, then $G$ contains a $K_{h,h,h}$-tiling.
%    \end{itemize}
% Furthermore, for $N$ sufficiently large and divisible by $h$, there % exists a $G_1\in{\cal G}_3(N)$
% such that $\bar{\delta}(G)=h\left\lceil 2N/(3h)\right\rceil+h-3$
% and $G$ has no $K_{h,h,h}$-tiling.
% In addition, for $N$ sufficiently large and divisible by $6h$,
% there exists a $G_0\in{\cal G}_3(N)$
% such that $\bar{\delta}(G)=h\left\lceil
% 2N/(3h)\right\rceil+h-2=(2/3)N+h-2$
% and $G$ has no $K_{h,h,h}$-tiling, which is best possible in
% that case.
 \label{thm:trithm}
\end{theorem}

So, the result is tight for $N=6h$, almost tight unless $N$ is an odd multiple of $6$ and, in the worst case, the upper and lower bounds
differ by $h$.

Clearly the complete tripartite graph $K_{h,h,h}$ can be perfectly
tiled by any 3-colorable graph on $h$ vertices. Since $f(h)\leq\frac{2N}{3}+2h-1$ whenever $N$ is divisible by $h$, we have the
following corollary.

\begin{corollary}\label{cor:main}
   Let $H$ be a $3$-colorable graph of order $h$. There
exists a positive integer $N_0$ such that if $N\geq N_0$ and $N$
divisible by $h$, then every $G\in {\cal G}_3(N)$ with
$\bar{\delta}(G)\ge \frac{2N}{3} + 2h -1$ contains an $H$-factor.
\end{corollary}

It is well known that every graph $G$ on $n=Nr$ vertices contains a
subgraph $G'\in {\cal G}_r(N)$ with $\bar{\delta}(G)\ge \delta(G)/r
- o(n)$ (following from a random balanced partition of the vertices
of $G$). Consequently Corollary~\ref{cor:main}  gives another proof
of the Alon-Yuster theorem \cite{AlYu2} for 3-colorable graphs as
follows: Fix a 3-colorable graph $H$ of order $h$ and let $G$ be a
graph of order $n=3N$ such that $N$ is sufficiently large and
divisible by $h$. If $\delta(G)\ge 2n/3 + \eps n$ for some $\eps>0$,
then we first find a subgraph $G'\in {\cal G}_3(N)$ with
$\bar{\delta}(G)\ge (2/3)N+ 2h-1$, and then apply
Corollary~\ref{cor:main} to $G'$ obtaining an $H$-factor in $G'$,
hence in $G$.

The proof of Theorem~\ref{thm:trithm} naturally falls into two parts as those of other tiling results
\cite{KSSz-AY, AliYi, MM, MSz}. In the first stage, we prove a
result that resembles the stability theorem of
Simonovits~\cite{Simo}; namely, any balanced tripartite graph with a
slightly weaker degree condition either contains an $K_{h,h,h}$-factor, or is
in a class of extremal graphs. In the second stage, we show that any
graph close to the extremal graphs contains an $K_{h,h,h}$-factor. This
approach seems to be a useful tool for obtaining exact results on
graphs or hypergraphs \cite{KSS-posa, KS1, KS2, MM, MSz}.  Our second stage turns out to be lengthy and intricate due to the fact that we must ensure that, when sets are partitioned, they must be divisible by $h$.

The structure of the paper is as follows. After stating the
Regularity Lemma and Blow-up Lemma in Section~\ref{sec:reglem}, we
prove the so-called ``fuzzy'' case (Theorem~\ref{thm:trifuzz}) in
Section~\ref{sec:fuzzy} and the extreme case
(Theorem~\ref{thm:triextr}) in Section~\ref{sec:extreme}.  The
graphs that provide the lower bounds for $f(h)$ in Theorem~\ref{thm:trithm} are constructed in
Section~\ref{sec:lb}.
% We give concluding remarks and open
% problems in Section~\ref{sec:sum}.

\section{Tools and Definitions}
\subsection{The Regularity Lemma and Blow-up Lemma}
\label{sec:reglem}

The Regularity Lemma and the Blow-up Lemma are main tools in the
proof of the so-called ``fuzzy'' case.  Let ``$+$'' designate a
disjoint union of sets.  We define the usual concepts of
$\epsilon$-regularity and $(\epsilon,\delta)$-super-regularity and
state the version of the Regularity Lemma that we use.
See~\cite{Reg-survey} and~\cite{Reg-survey2}.  In this paper, when
floors and ceilings are not crucial and do not effect the result, we
ignore them.

\begin{definition} The bipartite graph $G=(A,B,E)$ is {\bf
$\epsilon$-regular} if
\[ X\subset A, \qquad Y\subset B, \qquad |X|>\epsilon|A|,
   \qquad |Y|>\epsilon|B| \]
imply $|d(X,Y)-d(A,B)|<\epsilon$, otherwise we say $G$ is {\bf
$\epsilon$-irregular}. \label{def1}
\end{definition}

\begin{definition} $G=(A,B,E)$ is {\bf
$(\epsilon,\delta)$-super-regular} if
\[ X\subset A, \qquad Y\subset B, \qquad |X|>\epsilon|A|,
   \qquad |Y|>\epsilon|B| \]
imply $d(X,Y)>\delta$ and
\[ \deg(a)>\delta|B|,~\forall a\in A\qquad\mbox{and}\qquad
\deg(b)>\delta|A|,~\forall b\in B . \] \label{def2}
\end{definition}

\begin{lemma}[Regularity Lemma - Degree Form]
   For every positive $\epsilon$ there is an $M=M(\epsilon)$ such
   that if $G=(V,E)$ is any graph and $d\in [0,1]$ is any real
   number, then there is a partition of the vertex set $V$ into
   $\ell+1$ clusters $V_0,V_1,\ldots,V_{\ell}$ and there is a
   subgraph $G'=(V,E')$ with the following properties:
   \begin{itemize}
      \item $\ell\leq M$,
      \item $|V_0|\leq\epsilon|V|$,
      \item all clusters $V_i$, $i\geq 1$, are of the same size
            $L\leq\epsilon|V|$,
      \item $\deg_{G'}(v)>\deg_G(v)-(d+\epsilon)|V|$, $\forall
      v\in V$,
      \item $\left. G'\right|_{V_i}=\emptyset$ ($V_i$ are
      independent in $G'$),
      \item all pairs $\left. G'\right|_{V_i\times V_j}$,
      $1\leq i<j\leq l$, are $\epsilon$-regular, each with
      density either $0$ or exceeding $d$.
   \end{itemize}
\label{lem:reglem}
\end{lemma}

The proof of the regularity lemma (see~\cite{Sz}) begins with any
equipartition of the vertex set and refines it into a Szemer\'edi
partition, as defined above.  So, when we apply
Lemma~\ref{lem:reglem} to a balanced tripartite graph on $3N$
vertices with $N$ large enough, we can ensure that each cluster,
other than $V_0$, is a subset of exactly one piece of the
tripartition.

We will also need the Blow-up Lemma of Koml\'os, S\'ark\"ozy and
Szemer\'edi.  The graph $H$ can {\em be embedded into} graph $G$
if $G$ contains a subgraph isomorphic to $H$.
\begin{lemma}[Blow-up Lemma~\cite{KSSz-Blowup}]
   Given a graph $R$ of order $r$ and positive parameters $\delta$,
   $\Delta$, there exists an $\epsilon>0$ such that the following
   holds:  Let $N$ be an arbitrary positive integer, and let us
   replace the vertices of $R$ with pairwise disjoint $N$-sets
   $V_1,V_2,\ldots,V_r$ (blowing up).  We construct two graphs on
   the same vertex-set $V=\cup V_i$.  The graph $R(N)$ is obtained
   by replacing all edges of $R$ with copies of the complete
   bipartite graph $K_{N,N}$ and a sparser graph $G$ is constructed
   by replacing the edges of $R$ with some
   $(\epsilon,\delta)$-super-regular pairs.  If a graph $H$ with
   maximum degree $\Delta(H)\leq\Delta$ can be embedded into
   $R(N)$, then it can be embedded into $G$.
\label{lem:blowup}
\end{lemma}

\section{The Fuzzy Tripartite Theorem}
\label{sec:fuzzy}

The purpose of this section is to prove the following so-called
fuzzy tripartite theorem.
% Consider a tripartite graph $G$ with
% $\epsilon_0$ small and $N$ vertices, $N$ large enough, in each
% class satisfying
% \begin{equation}\label{eq:mindeg}
%   \bar{\delta}(G):= \min_{v\in V_i, j\neq i} \deg(v, V_j)\ge
%   \frac{2N}{3}
% \end{equation}
We say that $G$ is in {\em the extreme case with parameter $\Delta$}
if $G$ has three sets of size $\lfloor N/3\rfloor$, each in a
different vertex class, with pairwise density at most $\Delta$.
Recall that ${\cal G}_3(N)$ is the family of balanced tripartite
graphs with three parts, each of size $N$.

\begin{theorem}
   Given any positive integer $h$ and a $\Delta>0$, sufficiently
   small, there exists an $\epsilon>0$ and an integer
   $N_0=N_0(\Delta,h)$ such that whenever $N\geq N_0$, and $h$
   divides $N$, the following occurs: If
   $G=\left(V^{(1)},V^{(2)},V^{(3)};E\right)\in{\cal G}_3(N)$
   such that $\bar{\delta}(G)\geq
   (2/3-\epsilon)N$, then
   either $G$ has a subgraph which consists of $N/h$
   vertex-disjoint copies of the complete tripartite graph
   $K_{h,h,h}$ or $G$ is in the extreme case with parameter
   $\Delta$. \label{thm:trifuzz}
\end{theorem}

\subsubsubsection{Proof of Theorem~\ref{thm:trifuzz}}

As usual, there is a sequence of constants:
$$ \eps\ll\eps_1\ll \delta\ll d_1\ll\Delta_0\ll\Delta\ll h^{-1} . $$

% \[ \epsilon\ll\epsilon_1\ll\epsilon_5\ll
%    \epsilon_3\ll\alpha\ll\delta_4\ll\delta_3\ll d_3\ll d_1\ll
%    \epsilon_2\ll\Delta_0\ll\Delta \]

Begin with a tripartite graph
$G=\left(V^{(1)},V^{(2)},V^{(3)};E\right)$ with
$\left|V^{(1)}\right|=\left|V^{(2)}\right|=\left|V^{(3)}\right|=N$
in which each vertex is adjacent to at least $(2/3-\epsilon)N$
vertices in each of the other classes.
% Define the extreme case to be
% the case where $G$ has three sets of size $\lfloor N/3\rfloor$ with
% pairwise density at most $\Delta$.
Apply the Regularity Lemma (Lemma~\ref{lem:reglem}) with
$\epsilon_1$ and $d_1$, partitioning each $V_i$ into $\ell$ clusters
$V_1^{(i)}+\cdots+V_{\ell}^{(i)}$ of size $L\leq 3\epsilon_1 N$ and
an exceptional set $V_0^{(i)}$ of size at most $3\epsilon_1 N$. Let
us define $G_r$ to be the reduced graph defined as usual.
% It may be necessary to place some clusters into the
% exceptional sets to ensure that $\ell$ is divisible by 3.
It is important to observe that in the proof, the exceptional sets
will increase in size, but will always remain of size $O(\epsilon_1
N)$.

It is a routine calculation to see that the reduced graph $G_r$
(defined in the usual way where clusters are adjacent if the pair is
$\eps_1$-regular of density at least $d_1$) has the condition that
every cluster is adjacent to at least $(2/3-d_1)\ell$ clusters in
each of the other vertex classes.

\subsubsubsection{Step 1: Finding a triangle cover in $G_r$}

Here we can apply Lemma~\ref{lem:alcovlem}, the Almost-covering
Lemma, (Lemma 2.2 in~\cite{MM}) repeatedly to $G_r$ to get a
decomposition of $G_r$ into cluster-disjoint copies of $K_3$ and at
most $9$ clusters. If this is not possible, then $G$ must be in the
extreme case.

\begin{lemma}[Almost-covering Lemma~\cite{MM}]
   Let $H$ be a balanced tripartite graph on $3M$ vertices such
   that each vertex is adjacent to at least $(2/3-\epsilon)M$
   vertices in
   each of the other classes.  Let ${\cal T}_0$ be a partial $K_3$-tiling
   in $H$ with $|{\cal T}|<M-3$.  Then, either
   \begin{enumerate}
      \item there exists a partial $K_3$-tiling ${\cal T}'$ with
            $|{\cal T}'|>|{\cal T}|$ but
            $|{\cal T}'\setminus {\cal T}|\leq 15$, or
      \item $H$ has 3 subsets in 3 vertex classes of size
            $\lfloor M/3\rfloor$ with pairwise density at
            most $\Delta_0$.
   \end{enumerate}
\label{lem:alcovlem}
\end{lemma}

Note that the second case of Lemma~\ref{lem:alcovlem} implies that
$G_r$ is in the extreme case with parameter $\Delta_0$ and so $G$
itself is in the extreme case with parameter $\Delta\gg\Delta_0$.
The fact that $|{\cal T}'\setminus{\cal T}|\leq 15$ is not important
here but is crucial to arguments in Step 4.

We put the vertices from the clusters that are outside of the $K_3$-factor
into the corresponding exceptional set. For simplicity of notation,
we still denote the remaining graph by $G_r$ and assume that each
vertex class of $G_r$ has size $\ell$. The cluster-triangles which
cover $G_r$ are called $S_1, S_2, \ldots, S_{\ell}$, where
$S_i=\left\{S_i^{(1)},S_i^{(2)},S_i^{(3)}\right\}$ with
$S_i^{(j)}\subseteq V^{(j)}$ for $j=1,2,3$.

\subsubsubsection{Step 2: Making pairs in $S_i$ super-regular} For
each cluster-triangle, $S_i$, remove some vertices from it to make
each pair not just regular, but super-regular.  That is, remove a
vertex $v$ from a cluster in $S_i$ and place it in the exceptional
set if $v$ has fewer than $(d_1-\eps_1)L$ neighbors in each of the
other clusters of $S_i$. By $\eps_1$-regularity, there are at most
$2\eps_1 L$ such vertices in each cluster.  Remove more vertices if
necessary to ensure that each non-exceptional cluster is of the same
size, which is at least $(1-2\eps_1)L$ and divisible by $h$.

The Slicing Lemma is important for verifying that regularity is
maintained when small modifications are made to the clusters:
\begin{lemma}[Slicing Lemma, Fact 1.5 in~\cite{MM}]  Let $(A,B)$ be
an $\eps$-regular pair with density $d$, and, for some
$\alpha>\epsilon$, let $A'\subset A$, $|A'|\geq\alpha |A|$,
$B'\subset B$, $|B'|\geq\alpha |B|$.  Then $(A',B')$ is an
$\eps'$-regular pair with $\eps'=\max\{\eps/\alpha,2\eps\}$, and for
its density $d'$, we have $|d'-d|<\eps$.
\end{lemma}

Summarizing, any pair of clusters which was $\eps_1$-regular with
density at least $d_1$ is now $(2\eps_1)$-regular with density at
least $d_1-\eps_1$, as long as $\eps_1<1/4$. Furthermore, each pair
in a cluster-triangle $S_i$ is
$(2\eps_1,d_1-3\eps_1)$-super-regular. Each of the three exceptional
sets are now of size at most $\eps_1 N+\ell(2\eps_1 L)\leq 3\eps_1
N$.  The other clusters have the same number of vertices, which is
at least $(1-2\eps_1)L$ and is divisible by $h$.

\textbf{Remark:}  Because each triple
$\left(S_i^{(1)},S_i^{(2)},S_i^{(3)}\right)$, is super-regular,
we can apply the Blow-up Lemma to them (once we modify them to be of
equal size, divisible by $h$) so that they contain a spanning
subgraph of vertex-disjoint copies of $K_{h,h,h}$.

\subsubsubsection{Step 3: Create auxiliary triangles} In this step
we link each cluster to the corresponding one in the first
cluster-triangle, $S_1$. Its purpose is handling the last constant
number of leftover vertices in Step 5.

\begin{definition}
In a tripartite graph $G=\left(V^{(1)},V^{(2)},V^{(3)};E\right)$,
one vertex $x\in V^{(1)}$ (the cases of $x\in V^{(2)}$ or $V^{(3)}$
are defined similarly) is {\bf reachable} from another vertex $y\in
V^{(1)}$, if there is a chain of triangles $T_{1},\ldots,T_{2k}$
with $T_{j}=\left\{T_j^{(1)},T_j^{(2)},T_j^{(3)}\right\}$ and
$T_j^{(i)}\in V^{(i)}$ for $i=1,2,3$ such that the following occurs:
   \begin{enumerate}
      \item $x=T_1^{(1)}$ and $y=T_{2k}^{(1)}$,
      \item $T_{2j-1}^{(2)}=T_{2j}^{(2)}$ and
      $T_{2j-1}^{(3)}=T_{2j}^{(3)}$, for $j=1,\ldots,k$, and
      \item $T_{2j}^{(1)}=T_{2j+1}^{(1)}$, for $j=1,\ldots,k-1$.
   \end{enumerate}
\end{definition}

\begin{figure}
   \begin{center}
      \epsfig{file=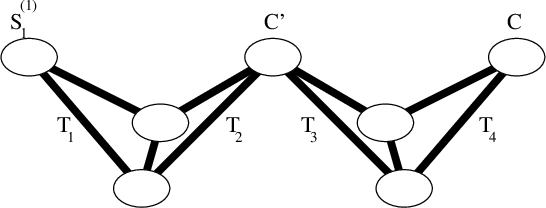}
   \end{center}
   \caption{An illustration of how cluster $S_1^{(1)}$ is reachable from a cluster $C$.}
   \label{fig:figSPAN}
\end{figure}

The Reachability Lemma (Lemma 2.6 in~\cite{MM}) states that, in the
reduced graph $G_r$, one cluster is reachable from any other cluster
in the same class using at most {\em four} cluster-triangles, unless
$G$ is in the extreme case. Thus, any cluster of $S_1$ is reachable
from every other cluster in the same class using at most 4
cluster-triangles (whose clusters come from at most $6$ different
$S_i$). Figure~\ref{fig:figSPAN} illustrates how $S_1^{(1)}$ is
reachable from cluster $C$.

Let $C$ be a cluster in $V^{(1)}$ and let $T_1,T_2$ or
$T_1,T_2,T_3,T_4$ be cluster-triangles that witness the fact that
$C$ is reachable from $S_1^{(1)}$.  Note that $T_1\cap
V^{(1)}=S_1^{(1)}$ and either both $k=1$ and $T_2\cap V^{(1)}=C$ or
$k=2$, $T_2\cap V^{(1)}=C'$ and $T_4\cap V^{(1)}=C$.

If $k=1$, then create a copy of $K_{h,h,h}$, called $H'$, in the
cluster triangle $T_1$, as well as an extra vertex in $C$ adjacent
to the vertices in $H'\cap V^{(2)}$ and $H'\cap V^{(3)}$.  This is
possible because all involved pairs are regular of nontrivial
density and $h$ is a constant compared with $L$, the size of
clusters.

If $k=2$, then create two copies of $K_{h,h,h}$.  The first is again
called $H'$, in the cluster triangle $T_1$, as well as an extra
vertex in the $V^{(1)}$ cluster that forms $T_2\cap T_3$, which is
adjacent to the vertices in $H'\cap V^{(2)}$ and $H'\cap V^{(3)}$.
The second copy of $K_{h,h,h}$, called $H''$, is in the cluster
triangle $T_3$ and there is a single vertex in $C$, which is
adjacent to the vertices in $H''\cap V^{(2)}$ and $H''\cap V^{(3)}$.

Color all of the vertices in $H'$ and in $H''$ (if it exists) {\em
red} and the additional vertex in $C$ and in $C'$ (if it exists)
{\em orange}.  If a vertex is not colored, we will heretofore call
it {\em uncolored}.  Repeat this $6h$ times for every cluster $C$ in
$G_r$, ensuring that all such red copies of $K_{h,h,h}$ and orange
vertices are vertex-disjoint.

This process of creating red copies of $K_{h,h,h}$ may result in a
discrepancy of uncolored vertices in the three clusters of some
$S_j$'s. A cluster may have at most $(3\ell)(6h)(h)=18\ell h^2$ red
vertices because there are $3\ell$ clusters $C$, the process is
iterated $6h$ times for each $C$ and no cluster gets more than $h$
vertices colored red with each iteration.  Similarly, no cluster
gets more than $\ell+1$ orange vertices.  We will remove some
uncolored vertices from each cluster, placing them in an exceptional
set.  This will be done to ensure that the number of uncolored
vertices in each cluster is the same and is divisible by $h$. Thus,
at most $18\ell h^2+(\ell+1)+(h-1)$ vertices are removed from any
cluster.  The sizes of exceptional sets are, thus, increased by at
most $(18\ell h^2+\ell+h)\ell\leq 20\ell^2 h^2$, a constant.

Summarizing, if we have $20\ell^2h^2\leq\eps_1 L$, then any pair
which, originally, was $\eps_1$-regular with density at least $d_1$
has that its uncolored vertices form a pair which is
$(4\eps_1)$-regular with density at least $d_1-2\eps_1$, as long as
$\eps_1<1/4$. Furthermore, the uncolored vertices in each pair of
some cluster-triangle $S_i$ is
$(2\eps_1,d_1-4\eps_1)$-super-regular. Each of the three exceptional
sets are now of size at most $3\eps_1 N+20\ell^2h^2\ell\leq 4\eps_1
N$.  The other (non-exceptional) clusters have at least
$(1-2\eps_1)L$ vertices with at most $\eps_1 L$ red vertices.  In
each non-exceptional cluster, the number of orange vertices is at
most $\eps_1 L$.  The number of uncolored vertices is at least
$(1-4\eps_1)L$ and is divisible by $h$.

\textbf{Remark:} This preprocessing ensures that we may later
transfer at most $6h$ vertices from any cluster $C$ to $S_1$ in the
following sense: Without loss of generality, suppose $C$ is a
cluster in $V^{(1)}$. In the case where $k=2$, identify an orange
vertex in $C$ and its corresponding red subgraphs $H'$ and $H''$ and
the orange vertex in the cluster $C'$. (The case where $k=1$ is
similar but simpler.)

Recolor the orange vertex in $C$ to be red.  Make a vertex from the
set $H''\cap C'$ to be uncolored and recolor the corresponding
orange vertex in $C'$ to be red.  Finally, uncolor a vertex from the
set $H'\cap S_1^{(1)}$. Except for $C$ and $S_1^{(1)}$, the number
of uncolored plus orange vertices remains the same in every cluster.
This number is decreased by 1 in $C$ but is increased by 1 in
$S_1^{(1)}$.  We will do this in Step 5.

\subsubsubsection{Step 4: Reducing the sizes of exceptional sets to
$6h$}
% Recall that the clusters in $G_r$ is covered by cluster-triangles
% $S_1,\ldots,S_{\ell}$. For $1\le i\le\ell$, we remove an uncolored vertex
% $v$ from a cluster in $S_i$ and place it in the exceptional set if
% $v$ has fewer than $\delta L$ (where $\delta=d_1-\eps_1$) neighbors
% among the uncolored vertices of either of the other two clusters in
% $S_i$. There are at most $2\eps_1 L$ such vertices in each cluster.
%  We may remove extra vertices from certain cluster such that each
% cluster has exactly $L'$ uncolored vertices where $L'$ is divisible by
% $h$. As a result, all three pairs in all $S_i$ become
% $(\eps_1,\delta/2)$ super-regular pairs and we may apply the Blow-up
% Lemma to tile $S_i$ with $K_{h,h,h}$ later.
Consider the exceptional sets $V_0^{(i)}$ for $i=1,2,3$, which are
all of the same size, at most $4\eps_1 N$.  We will show that we can
make them of size less than $6h$. So, suppose
$\left|V_0^{(1)}\right|\geq 6h$.

We will find $4h$ vertices in each exceptional set, bundling them into
$4$ sets of size $h$.  In the algorithm below, we will place at
least one bundle from each vertex class into $h$ vertex-disjoint
copies of $K_{h,h,h}$ and, together with at most $15$ additional
copies of $K_{h,h,h}$, we can remove them from the graph, reducing
the number of vertices in the exceptional set by at least $h$ vertices.

First, we observe that each vertex in the exceptional sets can be
regarded as a vertex in some non-exceptional cluster in the
following sense: Given an
$S_j=\left(S_j^{(1)},S_j^{(2)},S_j^{(3)}\right)$ and a vertex $v\in
V_0^{(1)}$. If $v$ is adjacent to at least $\delta L$ uncolored
vertices in $S_j^{(2)}$ and $S_j^{(3)}$, then we may remove a copy
of $K_{h,h,h}$ which consists of $v$, $h-1$ vertices from
$S_j^{(1)}$ and $h$ vertices from each of $S_j^{(2)}$ and
$S_j^{(3)}$.  This is easy to do because each pair is regular with
nontrivial density.
% The
% cluster $S_j^{(1)}$ thus have one more vertex than $S_j^{(2)}$ and
% $S_j^{(3)}$, which can be considered as a vertex in $V_0^{(1)}$.

Accordingly, we say a vertex $v\in V^{(i)}$ {\em belongs to} a
cluster $S_j^{(i)}$, for some, $i$ if $v$ is adjacent to at least
$\delta L$ uncolored vertices of each of the other clusters in
$S_j$. Using the minimum-degree condition, the number of clusters in
some other vertex class for which $v$ is adjacent to fewer than
$\delta L$ uncolored vertices is at most
\begin{equation}
   \frac{(1/3+\eps)N}{(1-3\epsilon_1)L-\delta L}\leq
   \frac{(1/3+\eps)\ell}{(1-3\epsilon_1-\delta)(1-\eps_1)} .
   \label{eq:delta}
\end{equation}
% We use $\eps\leq\eps_1\leq\delta/4$ and the fact that
% $(1-x)^{-1}\leq 1+3x/2$, if $x\leq 1/3$.
As long as $\delta$ is small enough, the expression in
(\ref{eq:delta}) is at most $(1/3+\delta)\ell$.  Thus, $v$ is
adjacent to at least $\delta L$ uncolored vertices in at least
$(2/3-\delta)\ell$ clusters in $V^{(k)}$, $k\neq i$. Hence, each
vertex in $V_0^{(i)}$ belongs to at least $(1/3-2\delta)\ell$
clusters.

Therefore, as long as $\left|V_0^{(1)}\right|\ge 3h$ and $\delta$ is
small enough ($\delta<1/(6h)$ is enough), the Pigeonhole Principle
guarantees that there is a cluster $C$ and set of $h$ vertices in
$V_0^{(1)}$ such that each these $h$ vertices belongs to $C$. Since
$\left|V_0^{(1)}\right|\geq 6h$, we can repeat this four times in
order to find four disjoint $h$-element subsets
$W_1^{(1)},\ldots,W_4^{(1)}$ of $V_0^{(1)}$ whose vertices belong to
disjoint clusters $C_1^{(1)},\ldots,C_4^{(1)}$, respectively.

Next we will show that, for $i=1,2,3$, there is some
$j\in\{1,2,3,4\}$ for which one of these sets $W_j^{(i)}$ of $h$
vertices can be {\em inserted} into $C_j^{(i)}$.

To see how this insertion works, it is useful to construct an
auxiliary graph $\tilde{G}$. The vertex set $V(\tilde{G})$ consists
of $\ell+4$ vertices in each of three vertex classes.  The first
$\ell$ vertices in each class represent the first $\ell$ clusters in
each class and two are adjacent if, originally, the pair of clusters
was $\eps_1$ regular with density at least $d_1$.  The remaining $4$
vertices in each class are merely duplicates of the vertices
representing the special clusters $C_1^{(i)},\ldots,C_4^{(i)}$ for
$i=1,2,3$.  Denote the duplicate of $C_j^{(i)}$ by
$\tilde{C}_j^{(i)}$.  They are adjacent to the same vertices as
their originals, but there are no edges between any of these $12$
duplicates.

Let $\mathcal{T}$ be the partial triangle-cover of $\tilde{G}$
corresponding to the triangle-cover $S_1,\ldots,S_{\ell}$ and apply
the Almost-covering Lemma (Lemma~\ref{lem:alcovlem}) to $\tilde{G}$
with this $\mathcal{T}$. The lemma provides a larger partial
triangle-cover $\mathcal{T}'$ which differs from $\mathcal{T}$ by at
most $15$ triangles.  We now create vertex-disjoint copies of
$K_{h,h,h}$ as follows: For each triangle in
$\mathcal{T}'\setminus\mathcal{T}$, find a copy of $K_{h,h,h}$ in
the uncolored vertices of the triple represented by that triangle.

To see how to deal with the duplicate clusters, suppose
$\tilde{C}_1^{(1)}$ is a vertex in
$\mathcal{T}'\setminus\mathcal{T}$. For each of the $h$ vertices
that belong to $\tilde{C}_1^{(1)}$, place it in a $K_{h,h,h}$ which
contains $h-1$ vertices in $\tilde{C}_1^{(1)}$ and $h$ in each of
the other clusters of the $S_i$ which contains $\tilde{C}_1^{(1)}$.
All of these copies of $K_{h,h,h}$ can be removed from the graph
entirely, they will be a part of the final $K_{h,h,h}$-factor of
$G$.  In the process of creating $\mathcal{T}'$, there may be a
cluster that was covered by $\mathcal{T}$ but is not covered by the
larger $\mathcal{T}'$.  In such a case, take an arbitrary set of $h$
uncolored vertices from that cluster and place it into the leftover
set.  Since $|\mathcal{T}'|>|\mathcal{T}|$, the net change in each
leftover set is the same and they each lose at least $h$ vertices.
Regardless, no cluster loses more than $h^2+h$ vertices.

We repeat this process until the number of vertices remaining in
each exceptional set is at most $6h$.  There is one caveat: If too
many vertices are removed from the clusters of $S_i$, then we will
not be able to apply the Blow-up Lemma later.  So, if in the process
of executing this algorithm, at least $(\delta/2)L$ vertices are
used from a cluster of $S_i$, then we say that $S_i$ is {\em dead}.

The number of dead cluster-triangles is not very large.  To see
this, there are three ways for vertices to leave a cluster.  First,
they are placed in a $K_{h,h,h}$ with a vertex from the leftover
set, so each vertex class $V^{(i)}$ loses at most $3|V_0^{(i)}|h$
vertices in this way. Second, each time $h$ vertices are inserted,
there are at most $15$ vertices that are a vertex in
$\mathcal{T}'\setminus\mathcal{T}$ and so could lose a total of
$15h$ vertices to a copy of $K_{h,h,h}$.  Third, there are at most
$3$ that are vertices uncovered by $\mathcal{T}'$ and so could lose
$3h$ vertices have to be placed from a cluster into the . Since this
algorithm is executed $|V_0^{(i)}|/h$ times, the total number of
vertices that leave clusters is at most
$$
3|V_0^{(i)}|h+\left(|V_0^{(i)}|/h\right)(15h+3h)=|V_0^{(i)}|(3h+18)\leq
6\epsilon_1 N(3h+18) . $$

The number of dead cluster triangles is at most
$$ \frac{6\epsilon_1 N(3h+18)}{(\delta/2) L}
=\frac{36(h+6)\epsilon_1}{\delta(1-\epsilon_1)}\ell . $$ So, as long
as $\epsilon_1\ll\delta\ll d_1$, the number of dead clusters is at
most $d_1\ell$ and the Almost-covering Lemma
(Lemma~\ref{lem:alcovlem}) can be applied to the live clusters
without changing the result because each cluster is adjacent to at
least $(2/3-d_1)\ell-d_1\ell$ clusters.

Summarizing, the graph $G$ consists of some vertex-disjoint copies
of $K_{h,h,h}$.  The remaining vertices induce a graph with clusters
that form triangles $S_1,\ldots,S_{\ell}$.  In each cluster, the
number of uncolored vertices is at least $(1-4\eps_1)L-\delta L\geq
(1-d_1/2)L$ and is divisible by $h$.
% The Slicing Lemma ensures that any pair that was
% $\eps_1$-regular with density at least $d_1$ has that
% its remaining uncolored vertices form an $(8\eps_1)$-regular
% pair with density at least $d_1-4\eps_1$.
For $i=1,\ldots,\ell$, the uncolored vertices in $S_i$ form a triple
that is pairwise $(8\eps_1,d_1-\delta)$-super-regular with density
at least $d_1-3\eps_1$, as long as $\delta\geq 4\eps_1$.  The edges
between other pairs of clusters are no longer relevant.

\subsubsubsection{Step 5: Inserting the last $\leq 6h$ leftover
vertices} Assume that each exceptional set has at most $6h$
vertices. Consider a vertex $x\in V_0^{(1)}$ and suppose $x$ belongs
to some cluster $C$. In the remark at the end of Step~2 we discussed
how to move a vertex from $C$ to $S_1^{(1)}$.  So, we place $x$ in
$C$, uncolor one of the orange vertices in $C$ and proceed by using
the red vertices and orange in the manner prescribed in Step~3,
until $S_1^{(1)}$ has an extra uncolored vertex.  Repeat until all
of the leftover vertices have been assigned to a cluster.

Then, uncolor all remaining orange vertices and remove the red
copies of $K_{h,h,h}$.  It remains to show that the vertices that
remain in each triple $S_i$ themselves form a
$K_{h,h,h}$-factor.  The number of uncolored vertices has not
decreased but has increased by at most $\eps_1 L$.  Now, we just
have to establish that uncolored vertices in the pairs of a
cluster-triangle remain super-regular for some set of parameters.

Each vertex in the leftover set is adjacent to at least
$(\delta/2)L$ vertices in the other live clusters if it belongs to
$C$ because it had been adjacent to $\delta L$ before Step~4 and at
most $(\delta/2)L$ vertices are removed by Step~4. Some
straightforward calculations, which we neglect to include here, show
that an $(8\eps_1)$-regular pair, with each cluster of size at least
$(1-d_1/2)L$, will be $(2\sqrt{\eps_1})$-regular if at most $\eps_1
L$ vertices are added to each set, as long as $\eps_1\ll d_1\ll
h^{-1}$.

Therefore, the pairs of vertices in each $S_i$ are
$(2\sqrt{\eps_1},\delta/2)$-super-regular and we can apply the
Blow-up Lemma to each $S_j$ to complete the $K_{h,h,h}$-factor of
$G$.

% Observe that for this case, it was not necessary to use the full
% degree condition, it was sufficient to have minimum degree at
% least $(2/3-\epsilon)N$.

\section{The Extremal Case}
\label{sec:extreme}
\newcommand{\Vexc}{The very extreme case}
\newcommand{\vexc}{the very extreme case}

Before we deal with the extremal case, we make the solution precise
by describing a specific exclusionary case, which we deal with in
Section~\ref{sec:vexc}.
\begin{definition}
A balanced tripartite graph $G$ on $3N$ vertices is in \vexc~if the
following occurs:  First, there are integers $N,q$ such that
$N=(6q+3)h$.  Second, there are sets $A^{(i)}_{j}\subseteq V^{(i)}$
for $i,j\in\{1,2,3\}$, each with size at least $2qh+1$, such that if
$v\in A^{(i)}_{j}$ then $v$ is nonadjacent to at most $3h-3$
vertices in $A^{(i')}_{j'}$ whenever the pair
$(A^{(i)}_{j},A^{(i')}_{j'})$ corresponds to an edge in the graph
$\Gamma_3$ with respect to the usual correspondence.
\end{definition}

The Main Theorem is proven by verifying the following:
\begin{theorem}
   Given any positive integer $h$, there exists a $\Delta$,
   $0<\Delta\ll h^{-1}$ and $N_0=N_0(h)$ such that whenever
   $N\geq N_0$ and $h$ divides $N$, the following occurs:
   If $G=\left(V^{(1)},V^{(2)},V^{(3)};E\right)$ is a balanced
   tripartite graph on $3N$ vertices and $G$ is in the extreme
   case with parameter $\Delta$ and $\bar{\delta}(G)\geq
   h\left\lceil \frac{2N}{3h}\right\rceil+h-1$, then, either $G$ has
   a $K_{h,h,h}$-factor or $N$ is an odd multiple of $3h$ and $G$ is in~\vexc.

   If $G$ is in~\vexc, we can find the $K_{h,h,h}$-factor if
   $\bar{\delta}(G)\geq h\left\lceil\frac{2N}{3h}\right\rceil+2h-1$. \label{thm:triextr}
\end{theorem}

Throughout all of  Section~\ref{sec:extreme}, assume that $G$ is
minimal, \ie, no edge of $G$ can be deleted so that the minimum
degree condition still holds. We will have the usual sequence of
constants:
$$ \Delta\ll\Delta_1\ll\Delta_2\ll
   \Delta_3\ll\Delta_4\ll\Delta_5\ll h^{-1} . $$
% \begin{eqnarray*}
%    \Delta & \ll & \delta_1\ll\delta_2\ll\delta_3\ll\delta_4
%    \ll\delta_5\ll\delta_6\ll\delta_7\ll\delta_8\ll\delta_9\ll\Delta_1
%    \ll\Delta_2\ll\Delta_3 \\
%    & & \ll\phi\ll\theta-7/8
% \end{eqnarray*}

We will assume for Parts 1, 2 and 3a
(Sections~\ref{sec:appfigONE},~\ref{sec:appfigTWO}
and~\ref{sec:theta}, respectively) that $\bar{\delta}(G)\geq
h\left\lceil\frac{2N}{3h}\right\rceil+h-1$.  In Part 3b
(Section~\ref{sec:gamma}), we will begin with the same assumption on
$\bar{\delta}$, until we are left with \vexc.  Then we will allow
$\bar{\delta}(G)\geq h\left\lceil\frac{2N}{3h}\right\rceil+2h-1$ in
Section~\ref{sec:vexc} to complete the proof.

\begin{definition}
   For $\delta$, $0<\delta<1$ a graph $H$ and positive integer
   $M$, we say a graph $G$ is {\em $\delta$-approximately $H(M)$}
   if $V(G)$ can be partitioned into $|V(H)|$ nearly-equally sized
   pieces, each corresponding to a vertex of $H$ so that for
   vertices $v,w\in V(H)$ with $v\not\sim_H w$, the parts of
   $V(G)$ corresponding to $v$ and $w$ have density between them
   less than $\delta$.
\end{definition}

\subsection{Part~1: The basic extremal case}
\label{sec:appfigONE}

For Part~1, we will prove that either a $K_{h,h,h}$-factor exists in
$G$, or $G$ is in Part~2.

Let $A^{(i)}\subset V^{(i)}$ for $i=1,2,3$ be the three pairwise
sparse sets given by the statement of the theorem and
$B^{(i)}=V^{(i)}\setminus A^{(i)}$ for $i=1,2,3$. We then define
$\widetilde{A}^{(i)}$ to be the ``typical'' vertices with respect to
$A^{(i)}$, $\widetilde{B}^{(i)}$ to be ``typical'' with respect to
$B^{(i)}$, and $C^{(i)}$ are what remain.  Formally, for $i=1,2,3$,
\begin{eqnarray*}
   \widetilde{A}^{(i)} & = & \left\{x\in V^{(i)} : \deg_{A^{(j)}}(x)\leq\Delta_1|A^{(j)}|,
   \forall j\neq i\right\} \\
   \widetilde{B}^{(i)} & = & \left\{y\in V^{(i)} :
   \deg_{A^{(j)}}(y)\geq(1-\Delta_1)|A^{(j)}|, \forall j\neq i\right\} \\
   C^{(i)} & = & V^{(i)}\setminus\left(A^{(i)}\cup B^{(i)}\right)
\end{eqnarray*}

As a result, we have that $|B^{(i)}\setminus
\widetilde{B}^{(i)}|\leq (2\Delta/\Delta_1)|B^{(i)}|$ and
$|A^{(i)}\setminus \widetilde{A}^{(i)}|\leq
(2\Delta/\Delta_1)|A^{(i)}|$.  So, with $\Delta_1=\Delta^{1/3}$,
$|\widetilde{B}^{(i)}|\geq \left(1-2\Delta_1^2\right)|B^{(i)}|\geq
\left(1-2\Delta_1^2\right)(2N/3)$ and $|\widetilde{A}^{(i)}|\geq
\left(1-2\Delta_1^2\right)\geq |A^{(i)}|\geq
\left(1-2\Delta_1^2\right)(N/3)$.  We ignore round-off in computing
sizes of $A^{(i)}$'s and $B^{(i)}$'s.

\subsubsubsection{Step~1: There are large $\widetilde{A}^{(i)}$
sets}

Let $t=h\left\lfloor N/(3h)\right\rfloor$.  We will eventually
modify each of the sets $\widetilde{A}^{(i)}$ into sets
$A^{(i)}_{1}$ that are either of size $t$ or $t+h$.  Let $N=(3q+r)h$
with $r\in\{0,1,2\}$.  More precisely, the largest $r$ sets
$\widetilde{A}^{(i)}$ will be modified into sets $A^{(i)}_{1}$ of
size $t+h$ and the smallest $3-r$ sets $\widetilde{A}^{(j)}$ will be
modified into sets $A^{(j)}_{1}$ of size $t$.

We will find, in
$\widetilde{A}^{(1)}\cup\widetilde{A}^{(2)}\cup\widetilde{A}^{(3)}$,
(vertex-disjoint) $h$-stars. We need the following lemma, proven in Section~\ref{sec:lemmas}.

\renewcommand{\labelenumi}{(\arabic{enumi})}
\begin{lemma}
   Let us be given $\epsilon>0$ and a positive integer $M$.
   \begin{enumerate}
      \item Let $(A^{(1)},A^{(2)};E)$ be a bipartite graph such that
      every vertex in $A^{(2)}$ is adjacent to at least $d_1$
      vertices in $A^{(1)}$.  Furthermore,
      $\left||A^{(i)}|-M\right|<\epsilon M$ and $d_i<\epsilon M$ for
      $i=1,2$.

      Provided $\epsilon<\left((h+1)h\right)^{-1}$, there is a
      family of vertex-disjoint copies of $K_{1,h}$ such that
      $\max\{0,d_1-h+1\}$ of them have centers in $A^{(1)}$.
      \label{lem:stars:bi}

      \item Let $(A^{(1)},A^{(2)},A^{(3)};E)$ be a tripartite graph such that
      every vertex not in $A^{(i)}$ is adjacent to at least $d_i$
      vertices in $A^{(i)}$, for $i=1,2,3$.  Furthermore,
      $\left||A^{(i)}|-M\right|<\epsilon M$ and $d_i<\epsilon M$ for
      $i=1,2,3$. \label{lem:stars:tri}

      Provided $\epsilon<\left(2(h+2)(h+1)h\right)^{-1}$, there
      is a family of vertex-disjoint copies of $K_{1,h}$ such
      that $\max\{0,d_i-h+1\}$ of them have centers in $A^{(i)}$ and
      leaves in $A^{(i+1)}$ (index arithmetic is modulo 3).
   \end{enumerate}
   \label{lem:stars}
\end{lemma}

With our degree condition, we can guarantee that each vertex not in
$V^{(i)}$ is adjacent to at least $|\widetilde{A}^{(i)}|-t+h-1$
vertices in $\widetilde{A}^{(i)}$. So, we use
Lemma~\ref{lem:stars}{\it (\ref{lem:stars:tri})} with $d_i\geq
|\widetilde{A}^{(i)}|-t+h-1$ to construct the stars with the
property that there are exactly enough centers in
$\widetilde{A}^{(i)}$ such that, when removed, the resulting set has
its size bounded above by either $t$ or $t+h$, whichever is
required. Place these centers into $Z^{(i)}$.

\subsubsubsection{Step~2: There are small $\widetilde{A}^{(i)}$
sets}

For a subgraph $K_{1,h,h}$, with $h\geq 2$, define the {\em
center} to be the vertex that is adjacent to all others. We will
also refer to the remaining vertices as {\em leaves}, although
their degree is $h+1$.

We will find, in $B:=\bigcup_{i=1}^3\left(\widetilde{B}^{(i)}\cup
C^{(i)}\right)$, (vertex-disjoint) copies of $K_{1,h,h}$ such that
$\max\{t+h-|\widetilde{A}^{(i)}|,0\}$ copies having its center
vertex in $B^{(i)}$ for the largest $r$ sets $\widetilde{A}^{(i)}$
and such that $t-|\widetilde{A}^{(j)}|$ copies having the center
vertex in $B^{(j)}$ for the smallest $3-r$ sets
$\widetilde{A}^{(j)}$. This will be accomplished with Lemma~\ref{lem:superstars}, proven in Section~\ref{sec:lemmas}.

\begin{lemma}
   Given $\delta>0$, there exists an
   $\epsilon=\epsilon(\delta)>0$ such that the
   following occurs:

   Let $(B^{(1)},B^{(2)},B^{(3)};E)$ be a tripartite graph such that for all
   $i\neq j$, each vertex in $B^{(i)}$ is adjacent to at least
   $(1-\epsilon)M$ vertices in $B^{(j)}$.  Furthermore,
   $\left||B^{(i)}|-2M\right|<\epsilon M$.

   If $(B^{(1)},B^{(2)},B^{(3)};E)$ contains no copy of $K_{1,h,h}$ with 1
   vertex in $B^{(1)}$, and $h$ vertices in each of $B^{(2)}$ and $B^{(3)}$,
   then the graph $(B^{(1)},B^{(2)},B^{(3)};E)$ is $\delta$-approximately
   $\Theta_{3\times 2}(M)$.
   \label{lem:superstars}
\end{lemma}

Lemma~\ref{lem:superstars} can be repeatedly applied at most
$\lceil\Delta_1(N/3)\rceil$ times, unless $G$ is
$\Delta_2$-approximately $\Theta_{3\times 3}(t)$. Here, we will want
$\Delta_1+6\Delta_1^2<\epsilon(\Delta_2)$.  Add the center vertices
of the $K_{1,h,h}$ subgraphs to the appropriate sets
$\widetilde{A}^{(i)}$.

Place vertices from $C^{(i)}$ into the sets $\widetilde{A}^{(i)}$ so
that $A^{(i)}_{1}$ is of size $t$ or $t+h$, for $i=1,2,3$ and that
$\sum_{i=1}^3\left|A^{(i)}_{1}\right|=N$.  Relabel the modified sets
$\widetilde{A}^{(i)}$ with $A^{(i)}_{1}$.

\subsubsubsection{Step~3: Finding a $K_{h,h}$-factor in $B$}

Now we try to find a $K_{h,h}$-factor among the remaining vertices
in $B$ with the goal of matching them with the $A^{(i)}_{1}$ vertices.
There are, however, some adjustments that should be made.

\begin{itemize}
   \item Vertices which are in copies of $K_{1,h,h}$, where
   the center vertex is in some $A^{(i)}_{1}$, will be in a
   specific copy of $K_{h,h}$ in $B$.
   \item If $v\in Z^{(i)}$ is the center of a $K_{1,h}$ with leaves
   in $A^{(i)}_{k}$, then $v$ will be assigned to $B^{(j)}$,
   where $\{j\}=\{1,2,3\}\setminus\{i,k\}$.
   \item Vertices $v\in C^{(i)}$ will be assigned to $B^{(j)}$
   if $v$ is adjacent to at least $(2\Delta_1)(N/3)$ vertices
   in $A^{(k)}$.  Since $v\in C^{(i)}$ it will be assigned either to
   $B^{(j)}$ or to $B^{(k)}$, where $\{j,k\}=\{1,2,3\}\setminus\{i\}$.
\end{itemize}

This last statement results from the fact that if $v\in C^{(i)}$, then
we may assume, without loss of generality, that $v$ is adjacent to
less than $(1-2\Delta_1^2)(2N/3)$ vertices in, say, $B^{(j)}$.  Hence,
$v$ is adjacent to at least $(2\Delta_1^2)(N/3)$ vertices in $A^{(j)}$
and at least $(3\Delta_1/2)(N/3)$ vertices in $A^{(i)}_{1}$.

Moreover, we have that $|C^{(i)}|\leq 9\Delta_1^2(N/3)$, $|Z^{(i)}|\leq
6\Delta_1^2(N/3)$ and there are at most $4\Delta_1^2(N/3)$ copies
of $K_{1,h,h}$ with the center vertex in a given $A^{(i)}_{1}$.

Lemma~\ref{lem:partition} is proven in Section~\ref{sec:lemmas}.

\begin{lemma}
   Let us be given $\delta>0$. Then there exists an
   $\epsilon=\epsilon(\delta)>0$ and a positive integer
   $t_0=t_0(\delta)$ such that the following occurs:

   Let there be positive integers $t_1,t_2,t_3$ which are divisible
   by $h$ and with $|t_i-t_j|\in\{0,h\}$, for all
   $i,j\in\{1,2,3\}$ and $t_1>t_0$.  Let $(B^{(1)},B^{(2)},B^{(3)};E)$ be a
   tripartite graph such that for distinct indices
   $i,j,k\in\{1,2,3\}$, $|B^{(i)}|=t_j+t_k$.  For all $i\neq j$,
   each vertex in $B^{(i)}$ is adjacent to at least
   $(1-\epsilon)t_1$ vertices in $B^{(j)}$.  We attempt to find
   a $K_{h,h}$-factor in the graph induced
   by $(B^{(1)},B^{(2)},B^{(3)};E)$ with certain restrictions:

   For each pair $(B^{(i)},B^{(j)})$, there are at most $\epsilon t_1$
   copies of $K_{h,h}$ which must be part of any factor.  For each
   $B^{(i)}$, there are at most $\epsilon t_1$ vertices with the
   following property: $v$ can only be in copies of $K_{h,h}$ in
   the pair $(B^{(i)},B^{(j)})$ and $v$ is adjacent to at least
   $(1-\epsilon)t_1$ vertices in $B^{(i)}$.

   If such a factor cannot be found, then, without loss of
   generality, the graph induced by $(B^{(1)},B^{(2)},B^{(3)};E)$ can be
   partitioned such that $B^{(i)}=B^{(i)}[1]+B^{(i)}[2]$, $|B^{(i)}[1]|=t_1$
   for $i=1,2,3$ and
   $d\left(B^{(j)}[1],B^{(2)}[1]\right)\leq\delta$ and
   $d\left(B^{(j)}[2],B^{(2)}[2]\right)\leq\delta$ for $j=1,3$.
   \label{lem:partition}
\end{lemma}

\begin{figure}
   \begin{center}
      \epsfig{file=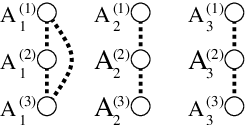}
   \end{center}
   \caption{The graph that defines Part 2.
   A dotted line represents a sparse pair.}
   \label{fig:figTWO}
\end{figure}

Then, match vertices in $C^{(i)}$ that are assigned to $B^{(j)}$ with $h$
typical neighbors in $B^{(j)}[i]$ and those with $h-1$ typical neighbors
in $B^{(i)}[j]$.  Finally, place the vertices that were moved into
copies of $K_{h,h,h}$.  All of these will be removed, allowing us to
apply Lemma~\ref{lem:partition}.  If the appropriate
$K_{h,h}$-factor cannot be found, then we are in the case of Part 2.
The diagram that defines that case is in Figure~\ref{fig:figTWO}.

\subsubsubsection{Step~4: Completing the $K_{h,h,h}$-factor}

We use Proposition~\ref{prop:sufffactor}, which allows us to
complete a $K_{h,h}$-factor into a $K_{h,h,h}$-factor.  The proof follows easily from K\"onig-Hall and is in Section~\ref{sec:lemmas}.

\begin{proposition} Let $h\geq 1$.
   \begin{enumerate}
      \item \label{it:bi:sufffactor}
      Let $G=(V^{(1)},V^{(2)};E)$ be a bipartite graph with
      $|V^{(1)}|=|V^{(2)}|=M$, $h$ divides $M$, and each vertex is
      adjacent to at least $\left(1-\frac{1}{2h^2}\right)M$
      vertices in the other part.  Then, we can find a $K_{h,h}$-factor in $G$.
      \item \label{it:tri:sufffactor}
      Let $G=(V^{(1)},V^{(2)},V^{(3)};E)$ be a tripartite graph
      with $|V^{(1)}|=|V^{(2)}|=|V^{(3)}|=M$, $h$ divides $M$, and each
      vertex is adjacent to at least
      $\left(1-\frac{1}{4h^2}\right)M$ vertices in
      each of the other parts.  Furthermore, let there
      be a $K_{h,h}$-factor in $(V^{(2)},V^{(3)})$.
      Then, we can extend it into a $K_{h,h,h}$-factor in $G$.
   \end{enumerate}
   \label{prop:sufffactor}
\end{proposition}

This allows us to find $K_{h,h,h}$-factors in each of
$\left(A^{(1)}_{1},B^{(2)}[3],B^{(3)}[2]\right)$,
$\left(A^{(2)}_{1},B^{(1)}[3],B^{(3)}[1]\right)$ and \linebreak
$\left(A^{(3)}_{1},B^{(1)}[2],B^{(2)}[1]\right)$ which completes the
$K_{h,h,h}$-factor in $G$.

\subsection{Part~2: $G$ is approximately the graph in
Figure~\ref{fig:figTWO}} \label{sec:appfigTWO}

{\bf Remark.}  In this part, we must deal with the fact that the
sets $A^{(2)}_2$ and $A^{(2)}_3$ may have close to the same number
of vertices, but that number is not divisible by $h$.  Much more
work needs to be done in order to modify these sets so that their
sizes become divisible by $h$.  We think it is easier to see the
basic arguments in the relatively shorter Part 1 before addressing
the specific issues raised in Part 2.

Recall that each vertex is adjacent to at least
$h\left\lceil\frac{2N}{3h}\right\rceil+h-1$ vertices in each of
the other pieces of the partition.  Again, let $t=h\lfloor
N/(3h)\rfloor$.  We will transform the graph that is
$\Delta_2$-approximately a graph defined by
Figure~\ref{fig:figTWO} with the vertices corresponding to sets of
size $N/3$.

Before we begin, we must examine the behavior of
$\left(A^{(1)}_{2}\cup A^{(1)}_{3},A^{(3)}_{2}\cup
A^{(3)}_{3}\right)$.  If this is $\Delta_5$- approximately
$\Theta_{2\times 2}(N/3)$, then call the dense pairs
$(E^{(1)},E^{(3)})$ and $(F^{(1)},F^{(3)})$.  Otherwise, coincidence
can only occur in either $V^{(1)}$ or $V^{(3)}$, but not both.
Without loss of generality, we will assume that if there is such a
coincidence, then it occurs in $V^{(1)}$.

We say that these pairs \textbf{coincide} with the sets
$A^{(i)}_{j}$ if the typical vertices of, say $A^{(3)}_{2}$, have
small intersection with those of $F^{(3)}$.  We will determine the
quantity that constitutes ``small'' later.  If $(E^{(1)},E^{(3)})$
and $(F^{(1)},F^{(3)})$ both coincide with
$(A^{(1)}_{2},A^{(3)}_{3})$ and $(A^{(1)}_{3},A^{(3)}_{2})$, then
$G$ is a graph that is approximately $\Theta_{3\times 3}(N/3)$
(Section~\ref{sec:theta}). If $(E^{(1)},E^{(3)})$ and
$(F^{(1)},F^{(3)})$ both coincide with $(A^{(1)}_{2},A^{(3)}_{2})$
and $(A^{(1)}_{3},A^{(3)}_{3})$, then approximately $\Gamma_3(N/3)$
(Section~\ref{sec:gamma}). Otherwise, coincidence can only occur in
either $V^{(1)}$ or $V^{(3)}$, but not both. Without loss of
generality, we will assume that if there is such a coincidence, then
it occurs in $V^{(1)}$.

Let $V^{(i)}=A^{(i)}_{1}+A^{(i)}_{2}+A^{(i)}_{3}+C^{(i)}$, such that
each $A^{(i)}_{j}$ has size $\left(1-3\Delta_2^{2/3}\right)t$ and
$\left(1+3\Delta_2^{2/3}\right)t$ and each vertex in $A^{(i)}_{j}$
is adjacent to at least $\theta t$ vertices in each set
$A^{(i')}_{j'}$ for which one of the following occurs:
\begin{itemize}
   \item $i=2$ and $j'\neq j$
   \item $i\in\{1,3\}$, $j=1$ and $j'\neq j$
   \item $i\in\{1,3\}$, $i'=2$ and $j\in\{2,3\}$
   \item $i\in\{1,3\}$, $i'=4-i$, $j\in\{2,3\}$ and $j'=1$
\end{itemize}
In other words, the vertices in $A^{(i)}_{j}$ are the ones that are
typical according to the rules established by
Figure~\ref{fig:figTWO}.  In addition, if, say $A^{(1)}_{2}$
coincides with $E^{(1)}$, then every vertex in $A^{(1)}_{2}$ is
adjacent to at least $\theta t$ vertices in $E^{(3)}$ and vice
versa. If there is no coincidence, then let $E^{(1)}$ and $E^{(3)}$
be redefined so that every vertex in $E_1$ is adjacent to at least
$\theta t$ vertices in $E^{(3)}$ and vice versa. Similarly for
$(F^{(1)},F^{(3)})$.

Each vertex $c\in C^{(2)}$ has the property that, for all
$j\in\{1,2,3\}$ and distinct $i',i''\in\{1,3\}$, if $c$ is
adjacent to fewer than $\Delta_3 t$ vertices in $A^{(i')}_{j}$, then
$c$ is adjacent to at least $\Delta_3 t$ vertices in $A^{(i'')}_{j}$.

Let $i\in\{1,3\}$, each vertex $c\in C^{(i)}$ has the property that,
for all $j\in\{1,2,3\}$,  $c$ cannot be adjacent to fewer than
$\Delta_3 t$ vertices in either $A^{(2)}_{2}$ or $A^{(2)}_{3}$.
Also, $c$ cannot be adjacent to fewer than $\Delta_3 t$ vertices in
both $A^{(2)}_{1}$ and $A^{(4-i)}_{1}$ or both $A^{(2)}_{2}$ and
$F^{(4-i)}$ (if it exists) or both $A^{(2)}_{3}$ and $E^{(4-i)}$ (if
it exists).

Trivially, each vertex in $V^{(i)}$ is adjacent to at least
$(1/2-\Delta_3)t$ vertices in at least two of
$\{A^{(i')}_{1},A^{(i')}_{2},A^{(i')}_{3}\}$ and in at least two of
$\{A^{(i'')}_{1},A^{(i'')}_{2},A^{(i'')}_{3}\}$, where $i',i''$ are distinct
members of $\{1,2,3\}\setminus\{i\}$.  This is particularly
important for vertices in $C^{(i)}$.

\subsubsubsection{Step 1: Ensuring small $A^{(i)}_{j}$ sets}

First, take each triple
$\left(A^{(1)}_{j},A^{(2)}_{j},A^{(3)}_{j}\right)$, $j=1,2,3$, and
construct disjoint copies of stars so that there are at most $t$
non-center vertices in each set $A^{(i)}_{j}$.  As in Part~1, we use
the fact that every vertex is adjacent to at least
$h\left\lceil\frac{2N}{3h}\right\rceil+h-1$ vertices in each of the
other parts as well as Lemma~\ref{lem:stars}.  For $i,j=1,2,3$,
place $|A^{(i)}_{j}|-t$ centers from $A^{(i)}_{j}$ into a set
$Z^{(i)}$.

\subsubsubsection{Step 2: Fixing the size of $A^{(i)}_{j}$ sets}

We have sets $A^{(i)}_{j}$ which have $|A^{(i)}_{j}|\leq t$ and the
remaining vertices are in sets $C^{(i)}\cup Z^{(i)}$.  Since $N$ is
divisible by $h$, we can place the vertices $C^{(i)}\cup Z^{(i)}$
arbitrarily into sets $A^{(i)}_{1}$, $A^{(i)}_{2}$ and $A^{(i)}_{3}$ so that
the resulting sets $A^{(i)}_{j}$ have cardinality $t$ or $t+h$ and for
$j=1,2,3$,
$$ |A^{(1)}_{j}|+|A^{(2)}_{j}|+|A^{(3)}_{j}|=N . $$
For this purpose, if $N/h\cong 1\pmod{3}$, add $h$ vertices to
each of $A^{(1)}_{2}$, $A^{(2)}_{1}$ and $A^{(3)}_{3}$.   If $N/h\cong
2\pmod{3}$, add $h$ vertices to all sets $A^{(i)}_{j}$, except
$A^{(1)}_{2}$, $A^{(2)}_{1}$ and $A^{(3)}_{3}$.

\subsubsubsection{Step 3: Partitioning the sets}

We will partition each set $A^{(i)}_{j}$ into two pieces, as close
as possible to equal size, but which have size divisible by $h$.
This must have  the property that a typical vertex in $A^{(i)}_{j}$
has at least $\left(1-2\Delta_4-6\Delta_2^{2/3}\right)(t/2)$
neighbors in each piece of the partition of $A^{(i')}_{j'}$, $i'\neq
i$, $j'\neq j$. Moreover, if a vertex has degree at least $\Delta_3
t$ in a set, it has degree at least $(\Delta_3/3)(t/2)$ in each of
the two partitions. Such a partition exists, almost surely, provided
$N$ is large enough, if the partition is random.

Assign to each part a permutation, $\sigma\in\Sigma_3$, which
assigns $j=\sigma(i)$. ($\Sigma_3$ denotes the symmetric group
that permutes the elements of $\{1,2,3\}$.)  Each part assigned to
$\sigma$ will be the same size.

% We call a vertex in $A^{(i)}_{j}$ a {\bf typical vertex} if it was not
% in $C^{(i)}$ and is neither a star-leaf nor a star-center.

% The equitable partition is verified in Appendix~\ref{app:soln}.

\subsubsubsection{Step 4: Assigning vertices}

The former $C^{(i)}$ vertices, as well as star-leaves and
star-centers, may only be able to form a $K_{h,h,h}$ with respect
to one particular permutation.

For example, consider a vertex $c$ which had been in $C^{(1)}$ but
is now in $A^{(1)}_{1}$.  Then, for either the pair
$(A^{(2)}_{2},A^{(3)}_{3})$ or the pair $(A^{(2)}_{3},A^{(3)}_{2})$,
the vertex $c$ is adjacent to at least $(1/2-\Delta_3)t$ in one set
and at least $\Delta_3 t$ vertices in the other; otherwise, it would
have been a typical vertex in $A^{(1)}_{1}$, $A^{(1)}_{2}$ or
$A^{(1)}_{3}$.

Assume that $c$ is adjacent to at least $\Delta_3 t$ vertices in
$A^{(2)}_{3}$ and at least $(1/2-\Delta_3)t$ vertices in $A^{(3)}_{2}$. In
this case, if $c$ were placed into the partition corresponding to
the identity permutation, then exchange $c$ with a typical vertex
in the partition assigned to $(23)$, using cycle notation of
permutations.

In a similar fashion, if there is a star with center in, say
$A^{(1)}_{2}$, and leaves in, say $A^{(2)}_{1}$, then we will use it to
form a $K_{h,h,h}$ with respect to the permutation $(12)\in
\Sigma_3$. Again, if any such leaf or center was in the wrong
partition, exchange it with a typical vertex in the other
partition.

The number of leaves in any set is at most
$2h\left(6\Delta_2^{2/3}t+h\right)$ and the number of centers is at
most $2\left(6\Delta_2^{2/3}t+h\right)$, the number of $C^{(i)}$
vertices is at most $9\Delta_2^{2/3}t$. So, if $N$ is large enough,
the total number of typical vertices in any $A^{(i)}_{j}$ which were
exchanged is at most $2(12h+21)\Delta_2^{2/3}t+4h^2+4h$.

With the partition established and the $C^{(i)}$, $Z^{(i)}$ and leaf
vertices in the proper part, we consider the triple formed by
three sets:
\begin{itemize}
   \item $A^{(2)}_{1}$, which will also be denoted $\widetilde{S}^{(2)}$
   \item the union of the piece of $A^{(1)}_{2}$ corresponding to
   $(12)$ and the piece of $A^{(1)}_{3}$ corresponding to $(132)$,
   denoted $\widetilde{S}^{(1)}$, and
   \item the union of the piece of $A^{(3)}_{2}$ corresponding to
   $(132)$ and the piece of $A^{(3)}_{3}$ corresponding to $(12)$,
   denoted $\widetilde{S}^{(3)}$.
\end{itemize}
Let the graph induced by the triple
$\left(\widetilde{S}^{(1)},\widetilde{S}^{(2)},\widetilde{S}^{(3)}\right)$
be denoted $\widetilde{\cal S}$.

\subsubsubsection{Step 5: Finding a $K_{h,h,h}$ cover in
$\widetilde{\cal S}$}

Let $t_0=|A^{(2)}_{1}|$.  First, take each $K_{1,h}$ in ${\cal S}'$
and complete it to form disjoint copies of $K_{h,h,h}$, using
unexchanged typical vertices.  This can be done if $\Delta_4$ is
small enough. Remove all such $K_{h,h,h}$'s containing stars.

Second, take each $c$ which had been a member of some $C^{(i)}$ and
use it to complete a $K_{h,h,h}$.  We can guarantee, because of the
random partitioning, that $c$ is adjacent to at least
$(\Delta_3/3)t_0$ vertices in one set and $(1/3-2\Delta_3)t_0$
vertices in the other.  Without loss of generality, let $c\in
\widetilde{S}^{(1)}$ with degree at least $(\Delta_3/3)t_0$ in
$\widetilde{S}^{(2)}$ and at least $(1/4-2\Delta_3)t_0$ in
$\widetilde{S}^{(3)}$. Since $\Delta_3\gg\Delta_2$, we can guarantee
$h$ neighbors of $c$ in $\widetilde{S}^{(2)}$ among unexchanged
typical vertices and, if $\Delta_3\ll\Delta_4\ll 1$, then $h$ common
neighbors of those among unexchanged typical vertices in
$N(c)\cap\widetilde{S}^{(3)}$. Finally, $\Delta_4\ll h^{-1}$ implies
this $K_{h,h}$ has at least $h-1$ more common neighbors in
$\widetilde{S}^{(1)}$. This is our $K_{h,h,h}$ and we can remove it.
Do this for all former members of a $C^{(i)}$.

Third, take each exchanged typical vertex and put it into a
$K_{h,h,h}$ and remove it.  Throughout this process, we have removed
at most $C_h\sqrt{\Delta_2}\times t_0$ vertices where $C_h$ is a
constant depending only on $h$.  What remains are three sets of the
same size, $t'\geq (1-C_h\sqrt{\Delta_2})t_0$, with each vertex in
$\widetilde{S}^{(1)}$ adjacent to at least, say
$\left(1/2-2\Delta_4\right)t'$, vertices in $\widetilde{S}^{(3)}$
and vice versa. Each vertex in $\widetilde{S}^{(1)}$ and in
$\widetilde{S}^{(3)}$ is adjacent to at least
$\left(1/2-2\Delta_4\right)t'$ vertices in $\widetilde{S}^{(2)}$ and
each vertex in $\widetilde{S}^{(2)}$ is adjacent to at least
$\left(1/2-2\Delta_4\right)t'$ vertices in $\widetilde{S}^{(1)}$ and
in $\widetilde{S}^{(3)}$.

Lemma~\ref{lem:Zhao}, from \cite{Zhao}, shows that we can find a
factor of $\left(\widetilde{S}^{(1)},\widetilde{S}^{(3)}\right)$
with vertex-disjoint copies of $K_{h,h}$ unless
$\left(\widetilde{S}^{(1)},\widetilde{S}^{(3)}\right)$ is
approximately $\Theta_{2\times 2}(N/6)$. In that case, find the
factor and finish to form a factor of $\widetilde{\cal S}$ of
vertex-disjoint copies of $K_{h,h,h}$ via K\"onig-Hall.
\begin{lemma}[Z.~\cite{Zhao}]
   For every $\epsilon>0$ and integer $h\geq 1$, there exists
   an $\alpha>0$ and an $N_0$ such that the following holds.
   Suppose that $N>N_0$ is divisible by $h$.  Then every
   bipartite graph $G=(A,B;E)$ with $|A|=|B|=N$ and
   $\delta(G)\geq (1/2-\alpha)N$ either contains a
   $K_{h,h}$-factor, or contains $A'\subseteq A$,
   $B'\subseteq B$ such that $|A'|=|B'|=N/2$ and
   $d(A,B)\leq \epsilon$.
   \label{lem:Zhao}
\end{lemma}

Lemma~\ref{lem:randpairs} states, in particular, that if a random
partition results in
$\left(\widetilde{S}^{(1)},\widetilde{S}^{(3)}\right)$ being
approximately $\Theta_{2\times 2}(N/6)$ with high probability, then
$\left(A^{(1)}_{2}\cup A^{(1)}_{3},A^{(3)}_{2}\cup
A^{(3)}_{3}\right)$ is approximately $\Theta_{2\times 2}(N/3)$.  The proof of Lemma~\ref{lem:randpairs} follows from similar arguments to those in the proof of Lemma 3.3 of~\cite{MM} and in Section 3.3.1 of~\cite{MSz} so we omit it.

\begin{lemma}
   For every $\epsilon>0$ and integer $h\geq 1$, there
   exists a $\beta>0$ and positive integer $t_0$ such that
   if $t\geq t_0$ the following holds.  Let $(A,B)$ be a
   bipartite graph such that
   $|A|,|B|\in\{2t-h,2t,2t+h\}$ with minimum degree at least
   $(1-\epsilon)t$ and is minimal with respect to this
   condition.  Let $A'\subset A$, $B'\subset B$,
   $|A'|=|B'|=t$ be chosen uniformly at random.  If
   $$ \Pr\{(A',B')\mbox{ contains a subpair with density
      at most }\epsilon\}\geq 1/4 $$
   then $(A,B)$ is $\beta$-approximately $\Theta_{2\times
   2}(t)$.
   \label{lem:randpairs}
\end{lemma}

We can, therefore, assume the existence of $(E^{(1)},E^{(3)})$ and
$(F^{(1)},F^{(3)})$.  Otherwise, Lemmas~\ref{lem:Zhao}
and~\ref{lem:randpairs} imply that $\widetilde{\cal S}$ has a
$K_{h,h,h}$-factor.

As a result, recall that we let the typical vertices in the dense
pairs in $\left(A^{(1)}_{2}\cup A^{(1)}_{3},A^{(3)}_{2}\cup
A^{(3)}_{3}\right)$ be denoted $(E^{(1)},E^{(3)})$ and
$(F^{(1)},F^{(3)})$. If the dense pairs do not coincide, then we
will work to ensure that
$|E^{(1)}\cap\widetilde{S}^{(1)}|=|E^{(3)}\cap S_3'|$ and
$|F^{(1)}\cap\widetilde{S}^{(1)}|=|F^{(3)}\cap\widetilde{S}^{(3)}|$
and both are divisible by $h$. Do this by moving vertices from
$\left(A^{(1)}_{2}\cap E^{(1)}\right)\setminus\widetilde{S}^{(1)}$
into $\left(A^{(1)}_{2}\cap E^{(1)}\right)\cap\widetilde{S}^{(1)}$
and move the same number from $\left(A^{(1)}_{2}\cap
F^{(1)}\right)\cap\widetilde{S}^{(1)}$ into $\left(A^{(1)}_{2}\cap
F^{(1)}\right)\setminus\widetilde{S}^{(1)}$. In addition, move
vertices from $\left(A^{(3)}_{2}\cap
E^{(3)}\right)\setminus\widetilde{S}^{(3)}$ into
$\left(A^{(3)}_{2}\cap E^{(3)}\right)\cap\widetilde{S}^{(3)}$ and
move the same number from $\left(A^{(3)}_{2}\cap
F^{(3)}\right)\cap\widetilde{S}^{(3)}$ into $\left(A^{(3)}_{2}\cap
F^{(3)}\right)\setminus\widetilde{S}^{(3)}$.

This can be done unless one of the intersections $A^{(i)}_{j}\cap
E^{(i)}$ or $A^{(i)}_{j}\cap F^{(i)}$ is too small.  This implies
the coincidence that we discussed at the beginning of this part. But
then, we have guaranteed that the remaining vertices of
$A^{(1)}_{2}$ are not only typical in that set but also typical in
$E^{(1)}$.  The same is true of $A^{(1)}_{3}$ and $F^{(1)}$.

Now, we want to move vertices in $V^{(3)}$ to ensure that
$|E^{(3)}\cap\widetilde{S}^{(3)}|=|A^{(1)}_{2}\cap\widetilde{S}^{(1)}|$
and
$|F^{(3)}\cap\widetilde{S}^{(3)}|=|A^{(1)}_{3}\cap\widetilde{S}^{(1)}|$.
Note that we have ensured that both
$|A^{(1)}_{2}\cap\widetilde{S}^{(1)}|$ and
$|A^{(1)}_{3}\cap\widetilde{S}^{(1)}|$ are divisible by $h$ and
approximately $N/6$.

We can do this as follows: Move vertices from $E^{(3)}\cap
A^{(3)}_{2}\setminus\widetilde{S}^{(3)}$ to $(E^{(3)}\cap
A^{(3)}_{2})\cap\widetilde{S}^{(3)}$ and move the same amount from
$(F^{(3)}\cap A^{(3)}_{2})\cap\widetilde{S}^{(3)}$ to $(F^{(3)}\cap
A^{(3)}_{2})\setminus\widetilde{S}^{(3)}$.  Also move vertices from
$(E^{(3)}\cap A^{(3)}_{3})\setminus\widetilde{S}^{(3)}$ to
$(E^{(3)}\cap A^{(3)}_{3})\cap\widetilde{S}^{(3)}$ and move the same
amount from $(F^{(3)}\cap A^{(3)}_{3})\cap\widetilde{S}^{(3)}$ to
$(F^{(3)}\cap A^{(3)}_{3})\setminus\widetilde{S}^{(3)}$.  Since none
of the intersections are small, this is possible. Complete this to
vertex-disjoint copies of $K_{h,h,h}$ in $\widetilde{\cal S}$ by
Proposition~\ref{prop:sufffactor}.

\subsubsubsection{Step 6: Completing the $K_{h,h,h}$-factor in $G$}

Now that we have found a $K_{h,h,h}$ that corresponds permutations
$(12)$ and $(132)$, we consider permutations in $\Sigma_3$. For a
$\sigma\in\Sigma_3\setminus\{(12),(132)\}$, let ${\cal
S}(\sigma)\stackrel{\rm
def}{=}\left(S^{(1)}_{\sigma(1)},S^{(2)}_{\sigma(2)},S^{(3)}_{\sigma(3)}\right)$
be a triple of parts formed by the random partitioning after the
exchange of vertices has taken place. The set $S^{(i)}_{\sigma(i)}$
is a subset of $A^{(i)}_{\sigma(i)}$. We have also ensured that
$s_{\sigma}\stackrel{\rm
def}{=}|S_{1,\sigma(1)}|=|S_{2,\sigma(2)}|=|S_{3,\sigma(3)}|$ and
$s_{\sigma}$ is divisible by $h$.  It is now easy to ensure that
this triple contains a $K_{h,h,h}$-factor:

First, take each star in ${\cal S}(\sigma)$ and complete it to
form disjoint copies of $K_{h,h,h}$, using unexchanged typical
vertices.  This can be done if $\Delta_4$ is small enough. Remove
all such $K_{h,h,h}$'s containing stars.

Second, take each $c$ which had been a member of some $C^{(i)}$ and
use it to complete a $K_{h,h,h}$.  We can guarantee, because of the
random partitioning, that $c$ is adjacent to at least
$(\Delta_3/3)s_{\sigma}$ vertices in one set and
$(2/3-2\Delta_3)s_{\sigma}$ vertices in the other.  Without loss of
generality, let $c\in S^{(1)}_{\sigma(1)}$ with degree at least
$(\Delta_3/3)s_{\sigma}$ in $S^{(2)}_{\sigma(2)}$ and at least
$(1/2-2\Delta_3)s_{\sigma}$ in $S^{(3)}_{\sigma(3)}$.  Since
$\Delta_3\gg\Delta_2$, we can guarantee $h$ neighbors of $c$ in
$S^{(2)}_{\sigma(2)}$ among unexchanged typical vertices and, if
$\Delta_3\ll\Delta_4\ll 1$, then $h$ common neighbors of those among
unexchanged typical vertices in $N(c)\cap S^{(3)}_{\sigma(3)}$.
Finally, $\Delta_4\ll h^{-1}$ implies this $K_{h,h}$ has at least
$h-1$ more common neighbors in $S^{(1)}_{\sigma(1)}$.  This is our
$K_{h,h,h}$ and we can remove it.  Do this for all former members of
a $C^{(i)}$.

Finally, take each exchanged typical vertex and put it into a
$K_{h,h,h}$ and remove it.  Throughout this process, we have
removed at most $C_h\sqrt{\Delta_2}\times s_{\sigma}$ vertices
where $C_h$ is a constant depending only on $h$.  What remains are
three sets of the same size, $s'\geq
(1-C_h\sqrt{\Delta_2})s_{\sigma}$, with each vertex adjacent to at
least, say $\left(1-2\Delta_4\right)s'$, vertices in each of the
other parts. If $N$ is large enough, then we can use the Blow-up
Lemma or
Proposition~\ref{prop:sufffactor}(\ref{it:tri:sufffactor}) to
complete the factor of ${\cal S}(\sigma)$ by copies of
$K_{h,h,h}$.

\subsection{Part~3a: $G$ is approximately $\Theta_{3\times
3}\left(\lfloor N/3\rfloor\right)$} \label{sec:theta}

\begin{figure}
   \begin{center}
      \epsfig{file=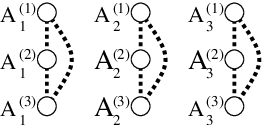}
   \end{center}
   \caption{The graph that defines Part 2.
   A dotted line represents a sparse pair.}
   \label{fig:figTH}
\end{figure}

Figure~\ref{fig:figTH} defines the case $\Theta_{3\times
3}(\lfloor N/3\rfloor)$ where sets that are connected with a
dotted line are sparse.

We will assume for this part that each vertex is adjacent to at
least $h\left\lceil\frac{2N}{3h}\right\rceil+h-1$ vertices in each
of the other pieces of the partition.  Again, let $t=h\lfloor
N/(3h)\rfloor$.

We will transform the $\Delta_2$-approximately $\Theta_{3\times
3}\left(\lfloor N/3\rfloor\right)$ by partitioning $V^{(i)}$,
$i=1,2,3$, into four sets, as follows:
$V^{(i)}=A^{(i)}_{1}+A^{(i)}_{2}+A^{(i)}_{3}+C^{(i)}$, such that each $A^{(i)}_{j}$ has
size between $(1-\sqrt{\Delta_2})t$ and $(1+\sqrt{\Delta_2})t$ and
each vertex in $A^{(i)}_{j}$ is adjacent to at least $\theta t$
vertices in each set $A^{(i')}_{j'}$ where $i'\neq i$ and $j'\neq j$.

Each vertex $c\in C^{(i)}$ has the property that, for all
$j\in\{1,2,3\}$ and distinct $i',i''\in\{1,2,3\}\setminus\{i\}$, if
$c$ is adjacent to fewer than $\Delta_3 t$ vertices in
$A^{(i')}_{j}$, then $c$ is adjacent to at least $\Delta_3 t$
vertices in $A^{(i'')}_{j}$; otherwise $c$ is in some set
$A^{(i)}_{j}$. Furthermore, $c$ is adjacent to at least
$(1/2-\Delta_3)t$ vertices in at least two of
$\left\{A^{(i')}_{1},A^{(i')}_{2},A^{(i')}_{3}\right\}$ and in at
least two of
$\left\{A^{(i'')}_{1},A^{(i'')}_{2},A^{(i'')}_{3}\right\}$.

\subsubsubsection{Step 1: Ensuring small $A^{(i)}_{j}$ sets}

First, take each triple
$\left(A^{(1)}_{j},A^{(2)}_{j},A^{(3)}_{j}\right)$, $j=1,2,3$, and
construct disjoint copies of stars so that there are at most $t$
non-center vertices in each set $A^{(i)}_{j}$. We use the fact that
every vertex is adjacent to at least
$h\left\lceil\frac{2N}{3h}\right\rceil+h-1$ vertices in each of the
other parts as well as Lemma~\ref{lem:stars}.  For $i,j=1,2,3$,
place $|A^{(i)}_{j}|-t$ centers from $A^{(i)}_{j}$ into a set
$Z^{(i)}$.

\subsubsubsection{Step 2: Fixing the size of $A^{(i)}_{j}$ sets}

We have sets $A^{(i)}_{j}$ which have $|A^{(i)}_{j}|\leq t$ and the
remaining vertices are in sets $C^{(i)}\cup Z^{(i)}$.  Since $N$ is
divisible by $h$, we can place the vertices $C^{(i)}\cup Z^{(i)}$
arbitrarily into sets $A^{(i)}_{1}$, $A^{(i)}_{2}$ and $A^{(i)}_{3}$ so that
the resulting sets $A^{(i)}_{j}$ have cardinality $t$ or $t+h$ and for
$j=1,2,3$,
$$ |A^{(1)}_{j}|+|A^{(2)}_{j}|+|A^{(3)}_{j}|=N . $$
For this purpose, we could place these vertices first to ensure
that all $|A^{(i)}_{j}|$ become of size exactly $t$.  If $N=3th+h$
then, for $i=1,2,3$, add all of the remaining $C^{(i)}\cup Z^{(i)}$ to
$A^{(i)}_{i}$.  If $N=3th+2h$ then, for $i=1,2,3$, add all of the
remaining $C^{(i)}\cup Z^{(i)}$ to $A^{(i)}_{j}$, $j\neq i$.

\subsubsubsection{Step 3: Partitioning the sets}

We will randomly partition each set $A^{(i)}_{j}$ into two pieces, as
close as possible to equal size but which have size divisible by
$h$, and assign them to a permutation, $\sigma\in \Sigma_3$, which
assigns $j=\sigma(i)$. ($\Sigma_3$ denotes the symmetric group
that permutes the elements of $\{1,2,3\}$.)  Each part assigned to
$\sigma$ will be the same size. We call a vertex in $A^{(i)}_{j}$ a
{\bf typical vertex} if it was not in $C^{(i)}$ and is neither a
star-leaf nor a star-center.

Note that a typical vertex in $A^{(i)}_{j}$ has at least
$(1-2\Delta_4-2\sqrt{\Delta_2})(t/2)$ neighbors in each piece of
the partition of $A^{(i')}_{j'}$, $i'\neq i$, $j'\neq j$, almost
surely -- provided $N$ is large enough and the partition was as
equitable as possible. Moreover, if a vertex has degree at least
$\Delta_3 t$ in a set, it has degree at least $(\Delta_3/3)(t/2)$
in each of the two partitions.

% The equitable partition is verified in Appendix~\ref{app:soln}.

\subsubsubsection{Step 4: Assigning vertices}

The former $C^{(i)}$ vertices, as well as star-leaves and star-centers
may only be able to form a $K_{h,h,h}$ with respect to one
particular permutation.

For example, consider a vertex $c$ which had been in $C^{(1)}$ but
is now in $A^{(1)}_{1}$.  Then, for either the pair
$(A^{(2)}_{2},A^{(3)}_{3})$ or the pair $(A^{(2)}_{3},A^{(3)}_{2})$,
the vertex $c$ is adjacent to at least $(1/2-\Delta_3)t$ in one set
and at least $\Delta_3 t$ vertices in the other.  It is easy to see
that, since $\Delta_2\ll \Delta_3$, that if this were not true, then
it would have been possible to place $c$ into one of the sets
$A^{(1)}_{1}$, $A^{(1)}_{2}$ or $A^{(1)}_{3}$.

Assume that $c$ is adjacent to at least $\Delta_3 t$ vertices in
$A^{(2)}_{3}$ and at least $(1/2-\Delta_3)t$ vertices in $A^{(3)}_{2}$. In
this case, if $c$ were placed into the partition corresponding to
the identity permutation, then exchange $c$ with a typical vertex
in the partition assigned to $(23)$, using cycle notation of
permutations.

In a similar fashion, if there is a star with center in, say
$A^{(1)}_{2}$, and leaves in, say $A^{(2)}_{1}$, then we will form a
$K_{h,h,h}$ with respect to the permutation $(12)\in \Sigma_3$.
Again, if any such leaf or center was in the wrong partition,
exchange it with a typical vertex in the other partition.

The number of leaves in any set is at most $2h(\sqrt{\Delta_2}
t+h)$ and the number of centers is at most $2(\sqrt{\Delta_2}
t+h)$, the number of $C^{(i)}$ vertices is at most $3\sqrt{\Delta_2}
t$. So, if $N$ is large enough, the total number of typical
vertices in any $A^{(i)}_{j}$ which were exchanged is at most
$(2h+6)\sqrt{\Delta_2} t$.

\subsubsubsection{Step~5: Completing the cover}

For some $\sigma\in \Sigma_3$, let ${\cal S}(\sigma)\stackrel{\rm
def}{=}\left(S^{(1)}_{\sigma(1)},S^{(2)}_{\sigma(2)},S^{(3)}_{\sigma(3)}\right)$
be a triple of parts formed by the random partitioning after the
exchange has taken place. The set $S^{(i)}_{\sigma(i)}$ is a subset
of $A^{(i)}_{\sigma(i)}$. We have also ensured in Step~3 that
$s_{\sigma}\stackrel{\rm
def}{=}|S^{(1)}_{\sigma(1)}|=|S^{(2)}_{\sigma(2)}|=|S^{(3)}_{\sigma(3)}|$
and $s_{\sigma}$ is divisible by $h$.  It is now easy to ensure that
this triple contains a $K_{h,h,h}$-factor:

First, take each star in ${\cal S}(\sigma)$ and complete it to
form disjoint copies of $K_{h,h,h}$, using unexchanged typical
vertices.  This can be done if $\Delta_4$ is small enough. Remove
all such $K_{h,h,h}$'s containing stars.

Second, take each $c$ which had been a member of some $C^{(i)}$ and
use it to complete a $K_{h,h,h}$.  We can guarantee, because of the
random partitioning, that $c$ is adjacent to at least
$(\Delta_3/3)s_{\sigma}$ vertices in one set and
$(2/3-2\Delta_3)s_{\sigma}$ vertices in the other.  Without loss of
generality, let $c\in S^{(1)}_{\sigma(1)}$ with degree at least
$(\Delta_3/3)s_{\sigma}$ in $S^{(2)}_{\sigma(2)}$ and at least
$(1/2-2\Delta_3)s_{\sigma}$ in $S^{(3)}_{\sigma(3)}$.  Since
$\Delta_3\gg\Delta_2$, we can guarantee $h$ neighbors of $c$ in
$S^{(2)}_{\sigma(2)}$ among unexchanged typical vertices and, since
$\Delta_3\ll\Delta_4\ll 1$, $h$ common neighbors of those among
unexchanged typical vertices in $N(c)\cap S^{(3)}_{\sigma(3)}$.
Finally, $\Delta_4\ll h^{-1}$ implies this $K_{h,h}$ has at least
$h-1$ more common neighbors in $S^{(1)}_{\sigma(1)}$.  This is our
$K_{h,h,h}$ and we can remove it.  Do this for all former members of
a $C^{(i)}$.

Finally, take each exchanged typical vertex and put it into a
$K_{h,h,h}$ and remove it.  Throughout this process, we have
removed at most $\Delta_2^{1/3}s_{\sigma}$ vertices if $\Delta_2$
is small enough. What remains are three sets of the same size,
$s'\geq \left(1-\Delta_2^{1/3}\right)s_{\sigma}$, with each vertex
adjacent to at least, say $\left(1-2\Delta_4\right)s'$, vertices
in each of the other parts. If $N$ is large enough, then we can
use the Blow-up Lemma or
Proposition~\ref{prop:sufffactor}(\ref{it:tri:sufffactor}) to
complete the factor of ${\cal S}(\sigma)$ by copies of
$K_{h,h,h}$.

\subsection{Part~3b: $G$ is approximately $\Gamma_3\left(\lfloor
N/3\rfloor\right)$} \label{sec:gamma}

\begin{figure}
   \begin{center}
      \epsfig{file=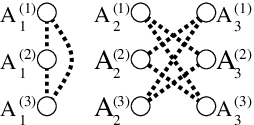}
   \end{center}
   \caption{The graph that defines Part 3b.
   A dotted line represents a sparse pair.}
   \label{fig:figGAM}
\end{figure}

Figure~\ref{fig:figGAM} defines the case $\Gamma_3(\lfloor
N/3\rfloor)$ where sets that are connected with a dotted line are
sparse.

We will assume for this part that each vertex is adjacent to at
least $h\left\lceil\frac{2N}{3h}\right\rceil+h-1$ vertices in each
of the other pieces of the partition.  We also assume that $G$ is
not in \vexc.  We must deal with \vexc~separately.

Let $t\stackrel{\rm def}{=}h\lfloor N/(3h)\rfloor$. We will
transform the $\Delta_2$-approximately $\Gamma_3\left(\lfloor
N/3\rfloor\right)$ by partitioning $V^{(i)}$, $i=1,2,3$, into four
sets, as follows: $V^{(i)}=A^{(i)}_{1}+A^{(i)}_{2}+A^{(i)}_{3}+C^{(i)}$, such that
each $A^{(i)}_{j}$ has size between $(1-\sqrt{\Delta_2})t$ and
$(1+\sqrt{\Delta_2})t$ and each vertex in $A^{(i)}_{1}$ is adjacent to
at least $(1-\Delta_3)t$ vertices in each set $A^{(i')}_{j'}$ where
$i'\neq i$ and $j'\in\{2,3\}$.  For $j=2,3$, $A^{(i)}_{j}$ is adjacent
to at least $(1-\Delta_3)t$ vertices in each set $A^{(i')}_{1}$ and
$A^{(i')}_{j}$, where $i'\neq i$.

Each vertex $c\in C^{(i)}$ has the property that, for all
$j\in\{1,2,3\}$ and distinct $i',i''\in\{1,2,3\}\setminus\{i\}$,
if $c$ is adjacent to fewer than $\Delta_3t$ vertices in
$A^{(i')}_{j}$, then $c$ is adjacent to at least $\Delta_3t$ vertices
in $A^{(i'')}_{j}$. Furthermore, $c$ is adjacent to at least
$(1/2-\Delta_4)t$ vertices in at least two of
$\left\{A^{(i')}_{1},A^{(i')}_{2},A^{(i')}_{3}\right\}$ and
$\left\{A^{(i'')}_{1},A^{(i'')}_{2},A^{(i'')}_{3}\right\}$.

Without loss of generality, we will assume that both
$|A^{(1)}_{2}|\geq |A^{(1)}_{3}|$ and $|A^{(2)}_{2}|\geq |A^{(2)}_{3}|$.

\subsubsubsection{Step 1: Ensuring small $A^{(i)}_{j}$ sets}

In each set $V^{(i)}$, we construct a set
$Z^{(i)}=Z^{(i)}[1]+Z^{(i)}[2]+Z^{(i)}[3]$ that will contain
star-centers.

If $|A^{(3)}_{2}|>|A^{(3)}_{3}|$, then $A^{(i)}_{2}$ is larger than
$A^{(i)}_3$ for $i=1,2,3$.  Use
Lemma~\ref{lem:stars}(\ref{lem:stars:bi}) to construct
$\max\left\{\min\{|A^{(i)}_{2}|-t\right., \left.t-|A^{(i)}_{3}|\},
0\right\}$ disjoint copies of $K_{1,h}$ in the
pair\footnote{Arithmetic in the indices is always done modulo 3.}
$(A^{(i)}_{2},A^{(i+1)}_{3})$ with centers in $A^{(i)}_{2}$. Place
these centers into $Z^{(i)}[3]$.

If $|A^{(3)}_{2}|<|A^{(3)}_{3}|$, we do something similar except
that first we use Lemma~\ref{lem:stars}(\ref{lem:stars:bi}) to
create the appropriate number of stars in
$(A^{(1)}_{2},A^{(2)}_{3})$ and $(A^{(2)}_{2},A^{(1)}_{3})$ with the
centers in $A^{(1)}_{2}$ and $A^{(2)}_{2}$, respectively. Place
these centers into $Z^{(1)}[3]$ and $Z^{(2)}[3]$, respectively.
Then, we apply Lemma~\ref{lem:stars}(\ref{lem:stars:bi}) to the pair
$(A^{(3)}_{3},A^{(2)}_{2})$.  (This $A^{(2)}_{2}$ is the possibly
modified set, with star-centers removed.)

By the conditions on Lemma~\ref{lem:stars}(\ref{lem:stars:bi}), we
see that each remaining set $A^{(i)}_{j}$ is of size at most $t$.
Now, apply Lemma~\ref{lem:stars}(\ref{lem:stars:tri}) to the triple
$\left(A^{(1)}_{1},A^{(2)}_{1},A^{(3)}_{1}\right)$.  For
star-centers in $A^{(i)}_{1}$, place $t-|A^{(i)}_{2}|$ into
$Z^{(i)}[2]$ and $t-|A^{(i)}_{3}|$ into $Z^{(i)}[3]$.

\subsubsubsection{Step 2: Fixing the size of the $A^{(i)}_{j}$ sets
for $j=1,2,3$}

\begin{figure}
\begin{center}
   \epsfig{file=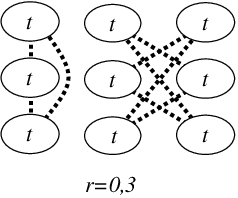}\hspace{1in}
   \epsfig{file=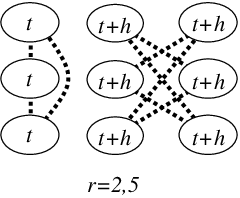}
\end{center}
   \caption{$N=(6q+r)h$ with $r=0,3$ and $r=2,5$, respectively.
   $t=2qh+h\lfloor r/3\rfloor$}
   \label{fig:figGAMt0325}
\end{figure}

\begin{figure}
\begin{center}
   \epsfig{file=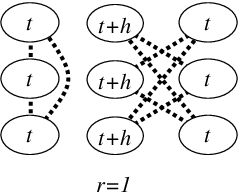}\hspace{1in}
   \epsfig{file=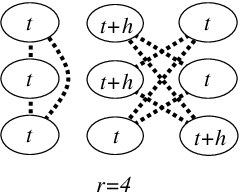}
\end{center}
   \caption{$N=(6q+r)h$ with $r=1$ and $r=4$, respectively.
   $t=2qh+h\lfloor r/3\rfloor$}
   \label{fig:figGAMt14}
\end{figure}

We now attempt to ``fill up'' the sets $A^{(i)}_{j}$.  Let $s_{i,j}$
be the targeted size.  There are several cases according to the
divisibility of $N/h$. Let $N/h=6q+r$ where $0\leq r<6$.
\begin{itemize}
   \item {\bf $r=0$:} $s_{i,j}=t$ for $i=1,2,3$ and $j=1,2,3$
   \item {\bf $r=1$:} $s_{i,j}=t$ for $i=1,2,3$ and $j=1,3$; and
   $s_{i,2}=t+h$ for $i=1,2,3$
   \item {\bf $r=2$:} $s_{i,1}=t$ for $i=1,2,3$; and
   $s_{i,j}=t+h$ for $i=1,2,3$ and $j=2,3$
   \item {\bf $r=3$:} $s_{i,j}=t$ for $i=1,2,3$ and $j=1,2,3$
   \item {\bf $r=4$:} $s_{i,1}=t$ for $i=1,2,3$; and
   $s_{1,3}=s_{2,3}=s_{3,2}=t$; and $s_{1,2}=s_{2,2}=s_{3,3}=t+h$
   \item {\bf $r=5$:} $s_{i,1}=t$ for $i=1,2,3$; and
   $s_{i,j}=t+h$ for $i=1,2,3$ and $j=2,3$
\end{itemize}
The cases of $r=0,3$ and $r=2,5$ are diagrammed in
Figure~\ref{fig:figGAMt0325} and the cases of $r=1$ and $r=4$ are
diagrammed in Figure~\ref{fig:figGAMt14}.

Place vertices of $Z^{(i)}[j]$ into $A^{(i)}_{j}$ for $i=1,2,3$ and
$j=1,2,3$. Furthermore, place vertices from $C^{(i)}$ into
$A^{(i)}_{j}$ for $i=1,2,3$ and $j=1,2,3$, ensuring that we still
have the case that $|A^{(i)}_{j}|\leq s_{i,j}$.

As usual, we call a vertex in $A^{(i)}_{j}$ a typical vertex if it was
neither in $C^{(i)}$ nor is either a star-leaf or a star-center. For
$j=2,3$, let ${\cal A}_j=\left(A^{(1)}_{j},A^{(2)}_{j},A^{(3)}_{j}\right)$. We
remove some copies of $K_{h,h,h}$ from among typical vertices of
these sets as follows:
\begin{itemize}
   \item {\bf $r=1$:} One from ${\cal A}_2$.
   \item {\bf $r=2$:} One from each of ${\cal A}_2$ and ${\cal
   A}_3$.
   \item {\bf $r=4$:} One from ${\cal A}_2$.
   \item {\bf $r=5$:} Two from ${\cal A}_2$.
\end{itemize}

Recalling $N=(6q+r)h$, each set is of size $2qh$ or $2qh+h$. Here
we note that $t_f\stackrel{\rm def}{=}h\lfloor t/(2h)\rfloor=qh$.
Also, $t_c\stackrel{\rm def}{=}h\lceil t/(2h)\rceil=qh$ if
$r=0,1,2$ and $t_c=(q+1)h$ if $r=3,4,5$.

\subsubsubsection{Step 3a: Partitioning the sets ($r\neq 3$)}

\newcommand{\tfloor}{t_f}
\newcommand{\tceil}{t_c}
Let $r\in\{0,1,2,4,5\}$.  Partition each $A^{(i)}_{1}$ set into parts
of nearly equal size.  Each part of the partition will receive a
label $\sigma\in\{1,2,3\}\times\{2,3\}$. Now, partition each
$A^{(i)}_{j}$ as follows:

\begin{quote}
Each $A^{(i)}_{1}$ will be split into two pieces: one of size
$\tfloor$ and another of size $\tceil$.  Unless both $r=4$ and
$i=3$, assign the smaller one with label $(i,2)$ and the larger
with label $(i,3)$.  If they are the same size, then assign them
arbitrarily.  If $r=4$ and $i=3$, then assign the one of size
$t_f$ with label $(3,3)$ and the one of size $t_c$ with $(3,2)$.

Each $A^{(i)}_{2}$ will be split into two pieces.  Unless both $r=4$
and $i\in\{1,2\}$, both pieces will be of size $\tfloor$ and will
be assigned $(i',2)$ and $(i'',3)$ arbitrarily, where
$\{i,i',i''\}=\{1,2,3\}$. If $r=4$ and $i\in\{1,2\}$, the one of
size $\tfloor$ is labeled $(3,2)$ and the one of size $\tceil$,
is labeled $(3-i,2)$.

Each $A^{(i)}_{3}$ will be split into two pieces.  Unless both $r=4$
and $i\in\{1,2\}$, both pieces will be of size $\tceil$ and will
be assigned $(i',2)$ and $(i'',3)$ arbitrarily, where
$\{i,i',i''\}=\{1,2,3\}$.  If $r=4$ and $i\in\{1,2\}$, the one of
size $\tfloor$ is labeled $(3,3)$ and one of size $\tceil$ is
labeled $(5-i,3)$.
\end{quote}
Figure~\ref{fig:figGAMtsplit} diagrams the partitioning.
\begin{figure}
\begin{center}
   \epsfig{file=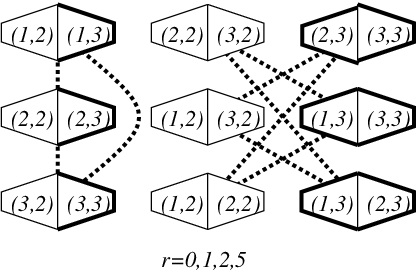}\hspace{.5in}
   \epsfig{file=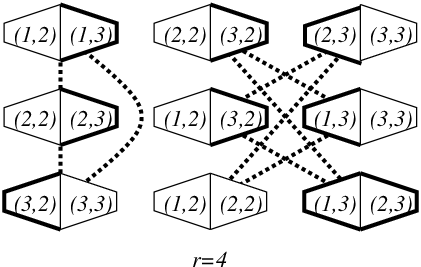}
\end{center}
   \caption{Partitioning the sets.  The light outlined half of a
   set is the piece of size $\tfloor$, the bold outlined half of a
   set is the piece of size $\tceil$.}
   \label{fig:figGAMtsplit}
\end{figure}

Partitioning the sets at random again ensures that the above can
be accomplished so that all of the vertices' neighborhoods
maintain roughly the same proportion, as in Part 3a, Step 3.

Now we proceed to Step 4.

\subsubsubsection{Step 3b: Partitioning the vertices ($r=3$, not
\vexc)}

Let $r=3$ (recall $N=(6q+r)h$) and let $G$ not be in \vexc. It may
be possible that there are additional stars $K_{1,h}$ between sparse
pairs. If it is possible to create enough such stars so as to move
star-centers into $Z^{(i)}$, then we can have at least one of these
sets $A^{(i)}_{j}$ of size at most $2qh$.  If we are not able to do
this, $G$ must be in the \vexc.  Without loss of generality, the set
to be made small is either $A^{(1)}_{1}$ or $A^{(1)}_{3}$.

\begin{itemize}
   \item Suppose vertices are removed to make $|A^{(1)}_{1}|=2qh$.
   We will make the set $A^{(1)}_{2}$ of size $(2q+2)h$ by adding
   vertices from the sets $C^{(1)}$, $Z^{(1)}[2]$ and $Z^{(1)}[1]$.

   \item Suppose vertices are removed to make $|A^{(1)}_{3}|=2qh$.
   We will make the set $A^{(1)}_{2}$ of size $(2q+2)h$ by adding
   vertices from the sets $C^{(1)}$, $Z^{(1)}[3]$ and $Z^{(1)}[1]$.
\end{itemize}

In each case, if the vertices in $Z^{(1)}[1]$ that were placed into
$A^{(1)}_{2}$ were themselves originally in $A^{(1)}_{2}$, then we
just treat them as typical vertices again, ignoring the star that
was formed.  Note that all sets are of size $(2q+1)h$, except
$|A^{(1)}_{2}|=(2q+2)h$ and either $A^{(1)}_{1}$ or $A^{(1)}_{3}$,
which has size $2qh$.  If $A^{(1)}_{1}$ is the small set, then
remove one copy of $K_{h,h,h}$ in the triple
$\left(A^{(1)}_{3},A^{(2)}_{3},A^{(3)}_{3}\right)$.

Now we partition each set as follows: Each $A^{(i)}_{1}$ will have
one piece of size $qh$ with label $(1,3)$. The other set will have
label $(1,2)$ size $(q+1)h$ in the case of $A^{(2)}_{1}$ and
$A^{(3)}_{1}$ and either $qh$ or $(q+1)h$ in the case of
$A^{(1)}_{1}$.  The set $A^{(1)}_{2}$ is partitioned into two pieces
of size $(q+1)h$, one labeled $(2,2)$, the other labeled $(3,2)$.
For $A^{(i)}_{2}$, $i=2,3$, we have one piece of size $qh$ and
labeled $(1,2)$ and the other of size $(q+1)h$, labeled $(5-i,2)$.
For $A^{(1)}_{3}$, it will have two pieces of size $qh$, one labeled
$(2,3)$, the other $(3,3)$. Finally, for $A^{(i)}_{3}$, $i=2,3$, we
have one piece of size $qh$ with label $(5-i,3)$ and the other will
have size either $qh$ or $(q+1)h$ and label $(1,3)$.

Partitioning the sets at random again ensures that the above can
be accomplished so that all of the vertices' neighborhoods
maintain roughly the same proportion, as in Part 3a, Step 3.

Now, we can proceed to Step 4.

\subsubsubsection{Step 4: Assigning vertices}

For any $\sigma\in\{1,2,3\}\times\{2,3\}$, we will show that the
$Z^{(i)}$ and $C^{(i)}$ vertices, in any $A^{(i)}_{j}$ can be assigned to one
of the two parts of the partition.

For example, consider a vertex $c$ which had been in $C^{(1)}$ but
is now in $A^{(1)}_{1}$.  Then, for either the pair
$(A^{(2)}_{2},A^{(3)}_{2})$ or the pair $(A^{(2)}_{3},A^{(3)}_{3})$,
the vertex $c$ is adjacent to at least $(1/2-\delta)t$ in one set
and at least $\Delta_3 t$ vertices in the other.  If such a pair is
$(A^{(2)}_{2},A^{(3)}_{2})$ then if $c$ were labeled $(1,2)$
exchange it with a typical vertex with label $(1,3)$.

Now, for example, consider a vertex $c$ which had been in $C^{(1)}$
but is now in $A^{(1)}_{2}$.  It is easy to check that for either
the pair $(A^{(2)}_{1},A^{(3)}_{2})$ or the pair
$(A^{(3)}_{1},A^{(2)}_{2})$, the vertex $c$ is adjacent to at least
$(1/2-\Delta_3)t$ in one set and at least $\Delta_3 t$ vertices in
the other.  If such a pair is, say, $(A^{(2)}_{1},A^{(3)}_{2})$, and
$c$ is not labeled $(2,2)$, then exchange it for a typical vertex of
that label.

Without loss of generality, this takes care of those vertices
$c\in C^{(i)}$.

Now we consider stars.  All star-centers are in sets $A^{(i)}_{2}$ or
$A^{(i)}_{3}$.  Without loss of generality, assume $z$ is such a
center in $A^{(1)}_{2}$ and the leaves are in $V^{(2)}$.  If the leaves
are in $A^{(2)}_{1}$, then $z$ must have been a member of $A^{(1)}_{1}$
originally.  So, $z$ and its leaves must have label $(2,2)$.  If
the leaves are in $A^{(2)}_{2}$, then $z$ must have been a member of
$A^{(1)}_{3}$ originally.  So, $z$ and its leaves must have label
$(3,2)$.  Exchange $z$ with typical vertices to ensure this.

Finally, we consider typical vertices moved from $A^{(i)}_{2}\cup
A^{(i)}_{3}$ to $A^{(i)}_{1}$.  Without loss of generality, suppose $z$ is
such a vertex in $A^{(1)}_{1}$.  If $z$ were originally from
$A^{(1)}_{2}$, then it is a typical vertex with respect to $A^{(2)}_{2}$
and $A^{(3)}_{2}$ and $z$ should receive label $(1,2)$.  Otherwise, it
is typical with respect to $A^{(2)}_{3}$ and $A^{(3)}_{3}$ and $z$ should
receive label $(1,3)$.

This completes the verification that all moved vertices can
receive at least one label of the $A^{(i)}_{j}$ set in which it is
placed.

\subsubsubsection{Step 5: Completing the cover}

For any $\sigma\in\{1,2,3\}\times\{2,3\}$, let ${\cal S}(\sigma)$
be one of the triples defined above.  We can finish as in Part 3a,
Step 5.

\subsection{\Vexc}
\label{sec:vexc}

Recall \vexc:
\begin{quote}
   There are integers $N,q$ such that $N=(6q+3)h$.  There are
   sets $A^{(i)}_{j}$ for $i,j\in\{1,2,3\}$, with sizes at least
   $2qh+1$, such that if $v\in A^{(i)}_{j}$ then $v$ is nonadjacent
   to at most $3h-3$ vertices in $A^{(i')}_{j'}$ whenever the pair
   $(A^{(i)}_{j},A^{(i')}_{j'})$ corresponds to an edge in the graph
   $\Gamma_3$ with respect to the usual correspondence.
\end{quote}

In this case, we must raise the minimum degree condition to
$2N/3+2h-1$.  Recalling Part 4, Step 3b, we were able to proceed if
we were able to make one of the sets $A^{(i)}_{j}$ small by means of
creating stars.  Each vertex in $A^{(2)}_{2}$ is adjacent to at
least $|A_{1,3}|-N/3+2h-1$ vertices in $A^{(1)}_{3}$.  Using
Lemma~\ref{lem:stars}(\ref{lem:stars:bi}), we have that there is a
family of $|A^{(1)}_{3}|-N/3+h$ vertex-disjoint stars with centers
in $A^{(1)}_{3}$.  We move the centers to $A^{(1)}_{2}$.  Then we
can proceed from Part 4, Step 4.

\subsection{Proofs of Lemmas}
\label{sec:lemmas}

\begin{proofcite}{Lemma~\ref{lem:stars}}
   \begin{enumerate}
      \item Let $\delta_1=d_1-h+1$.  If the stars cannot be
      created greedily, then
      there is a set $S\subset A^{(1)}$ and $T\subset A^{(2)}$ such
      that $|S|\leq\delta_1-1$ and $|T|=|S|h$ and each
      vertex in $A^{(1)}\setminus S$ is adjacent to less than $h-1$
      vertices in $A^{(2)}\setminus T$.  In this case,
      $$ (d_1-|S|)|A^{(2)}\setminus T|\leq e(A^{(1)}\setminus S,
      A^{(2)}\setminus T)\leq (h-1)|A^{(1)}\setminus S| . $$

      This gives
      \begin{eqnarray*}
         |S| & \geq & \delta_1
                      -(h-1)\frac{|A^{(1)}\setminus S|
                                  -|A^{(2)}\setminus T|}
                                 {|A^{(2)}\setminus T|} \\
               & \geq & \delta_1-(h-1)\frac{|A^{(1)}|-|A^{(2)}|+(h-1)|S|}
                                           {|A^{(2)}|-h|S|} \\
               & \geq & \delta_1-(h-1)\frac{(h+1)\epsilon M}
                                           {(1-(h+1)\epsilon)M} .
      \end{eqnarray*}
      If $\epsilon<(h^2+h)^{-1}$, then this gives
      $|S|>\delta_1-1$.  Since $|S|$ is an integer,
      $|S|\geq\delta_1$, contradicting the condition we put on
      $|S|$.

      \item Let $\delta_i=\max\{0,d_i-h+1\}$ for $i=1,2,3$.  If,
      say, $\delta_3=0$, then apply part (\ref{lem:stars:bi}) to
      the pair $(A^{(2)},A^{(3)})$ to create $\delta_2$ vertex-disjoint
      stars with centers in $A^{(2)}$.  Let $Z^{(2)}$ be the set of the
      centers.  Apply part (\ref{lem:stars:bi}) to
      $(A^{(1)},A^{(2)}\setminus Z^{(2)})$ and we can find $\delta_1$
      vertex-disjoint stars with centers in $A^{(1)}$ if
      $2\epsilon<(h^2+h)^{-1}$.

      So, we may assume that $\delta_i>0$ for $i=1,2,3$.  Note that
      if it is possible to construct $\delta_1+\delta_2$ disjoint
      copies of $K_{1,h}$ in $(A^{(1)},A^{(2)})$ with centers,
      $Z^{(1)}\subset A^{(1)}$, then we can finish with application of part
      (\ref{lem:stars:bi}).  To see this, apply part
      (\ref{lem:stars:bi}) to $(A^{(3)},A^{(1)}\setminus Z^{(1)})$,
      with
      $3\epsilon<(h^2+h)^{-1}$, creating $\delta_3$ stars with
      centers $Z^{(3)}\in A^{(3)}$.  Then apply part (\ref{lem:stars:bi})
      to $(A^{(2)},A^{(3)}\setminus Z^{(3)})$ ($2\epsilon<(h^2+h)^{-1}$).
      There will be $\delta_1$ stars remaining in $(A^{(1)},A^{(2)})$ which
      are vertex-disjoint from the rest.

      So, we will assume that it is not possible to create
      $\delta_1+\delta_2$ vertex-disjoint copies of $K_{1,h}$
      in $(A^{(1)},A^{(2)})$ with centers in $A^{(1)}$.  That means there is an
      $S\subset A^{(1)}$ and a $T\subset A^{(2)}$ such that
      $|S|<\delta_1+\delta_2$, $|T|=h|S|$ and every vertex in
      $A^{(1)}\setminus S$ is adjacent to at most $h-1$ vertices in
      $A^{(2)}\setminus T$.

      Now apply part (\ref{lem:stars:bi}) to $(A^{(3)},A^{(1)}\setminus S)$
      to obtain $\delta_3$ vertex-disjoint copies of $K_{1,h}$ with
      centers $Z^{(3)}\subset A^{(3)}$.  (Here, we need
      $3\epsilon<(h^2+h)^{-1}$.) Next, apply part
      (\ref{lem:stars:bi}) to $(A^{(2)},A^{(3)}\setminus Z^{(3)})$ to obtain
      $\delta_2$ vertex-disjoint copies of $K_{1,h}$ with centers
      $Z^{(2)}\setminus A^{(2)}$.  (Here, we need
      $2\epsilon<(h^2+h)^{-1}$.)  Finally, apply part
      (\ref{lem:stars:bi}) to $\left(A^{(1)},A^{(2)}\setminus (Z^{(2)}\cup
      T)\right)$ to obtain $\delta_1$ vertex-disjoint copies of
      $K_{1,h}$ with centers $Z^{(1)}\subset A^{(1)}$.  (Here, we need
      $(2h+2)\epsilon<(h^2+h)^{-1}$.)  But, because no vertex in
      $A^{(1)}\setminus S$ is adjacent to $h$ vertices in
      $A^{(2)}\setminus (Z^{(2)}\cup T)$, it must be the case that
      $Z^{(1)}\subset S$ and our $\delta_1+\delta_2+\delta_3$ copies of
      $K_{1,h}$ are, indeed, vertex-disjoint.
   \end{enumerate}
\end{proofcite}

\begin{proofcite}{Lemma~\ref{lem:superstars}}
   We can first apply the following theorem of
   Erd\H{o}s, Frankl and R\"odl~\cite{EFR}:
   \begin{theorem}
      For every $\epsilon'>0$ and graph $F$, there is a constant
      $n_0$ such that for any graph $G$ of order $n\geq n_0$, if
      $G$ does not contain $F$ as a subgraph, then $G$ contains a
      set $E'$ of at most $\epsilon' n^2$ edges such that
      $G\setminus E'$ contains no $K_r$ with $r=\chi(F)$.
      \label{thm:erdos}
   \end{theorem}
   Here, $F=K_{1,h,h}$ and $r=3$.

   So, after removing at most $\epsilon' (3M)^2$ edges, we have
   that the number of vertices in each part that are adjacent
   to at least $\sqrt{\epsilon}M$ vertices in each of the other
   two parts is at least
   $\left(1-\frac{18\epsilon'}{\sqrt{\epsilon}-\epsilon}\right)M$.

   So, now
   $2\left(1-\frac{9\epsilon'}{\sqrt{\epsilon}-\epsilon}\right)M\leq
   |B^{(i)}|\leq 2\left(1+\frac{\epsilon}{2}\right)M$ and each vertex
   is adjacent to at least
   $\left(1-\frac{18\epsilon'}{\sqrt{\epsilon}-\epsilon}\right)M$
   vertices in each of the other two parts.

   Finally, we use a version of a proposition appearing
   in~\cite{MM}, rephrased below:
   \begin{proposition}
      For a $\Delta$ small enough, there exists $\epsilon''>0$
      such that if $H$ is a tripartite graph with at least
      $2\left(1-\epsilon''\right)t$ vertices in each vertex class
      and each vertex is nonadjacent to at most
      $\left(1+\epsilon''\right)t$ vertices in each of the other
      classes. Furthermore, let $H$ contain no triangles.  Then,
      each vertex class is of size at most
      $2\left(1+\epsilon''\right)t$ and $H$ is $\Delta$-approximately
      $\Theta_{3\times 2}(t)$. \label{pTHETA}
   \end{proposition}

  By guaranteeing $\epsilon''\gg\epsilon'\gg\epsilon$ and
  $\delta=\Delta(\epsilon'')+\epsilon''$, the lemma follows.
\end{proofcite}

\begin{proofcite}{Lemma~\ref{lem:partition}}

   Let $\epsilon'$ be chosen such that $\epsilon'\ll\delta$.

   For this lemma, we partition the possibilities according to
   whether the pairs $(B^{(i)},B^{(j)})$ are approximately
   $\Theta_{2\times 2}(t_1)$.  That is, there are two sets of
   size $t_1$ which have density less than $\epsilon'$.
   Minimality gives the rest.

   In addition, we say that graphs $\Theta_{2\times 2}(t_1)$
   {\bf coincide} if there are sets $\widetilde{B}^{(i)}\subseteq B^{(i)}$,
   $\widetilde{B}^{(j)}\subseteq B^{(j)}$, $\widetilde{B}^{(k)}\subseteq
   B^{(k)}$, all of size
   $t_1$, such that both $(\widetilde{B}^{(i)},\widetilde{B}^{(j)})$ and $(\widetilde{B}^{(j)},\widetilde{B}^{(k)})$ have density less than $\epsilon'$.

   \subsubsubsection{Case 1: No pair is $\Theta_{2\times 2}(t_1)$}

   For each distinct $i,j,k\in\{1,2,3\}$, partition $B^{(i)}$ into
   two pieces, $B^{(i)}[j]$ and $B^{(i)}[k]$ with $|B^{(i)}[j]|=t_j$
   and $|B^{(i)}[k]|=t_k$.  If this partition is done uniformly at
   random, then with probability approaching 1, each vertex in
   $B^{(i)}[k]$ is adjacent to at least $(1/2-\epsilon^{1/2})t_k$
   vertices in $B^{(j)}[k]$.  So there exists a partition such
   that each vertex in $B^{(i)}$ is adjacent to at least
   $(1/2-\epsilon^{1/2})t_1$ vertices in each of the pieces
   $B^{(j)}[k]$, $j,k\neq i$ and such that the pair
   $(B^{(2)}[1],B^{(3)}[1])$ fails to contain a subpair with
   $\lfloor t_1/2\rfloor$ vertices in each part and density at
   most $\epsilon^{1/3}$.

   The vertices that are reserved will have to be placed in the
   proper set.  For example, if a reserved $K_{h,h}$ is in the
   pair $(B^{(i)},B^{(j)})$, then those vertices will need to be in the
   pair $(B^{(i)}[k],B^{(j)}[k]$.  So, we exchange vertices in
   $B^{(i)}[k]$ for vertices in $B^{(i)}[j]$ so that reserved
   vertices are in the proper place.  At most
   $4(\epsilon+\epsilon)t_1$ vertices are either reserved or
   moved in each set $B^{(i)}[j]$.  After such exchanges occur,
   place the moved vertices into vertex-disjoint copies of $K_{h,h}$
   that lie entirely within the given pairs.  This can be done
   because each vertex not in $B^{(i)}$ is adjacent to almost half
   of the vertices in both $B^{(i)}[j]$ and $B^{(i)}[k]$.

   Consider what remains of these sets.  The number of vertices
   is still divisible by $h$ and at most
   $8h(\epsilon)t_1$ have been placed into these
   copies of $K_{h,h}$.  We look for a perfect $K_{h,h}$-factor
   in each of the pairs $(B^{(1)}[3],B^{(2)}[3])$,
   $(B^{(1)}[2],B^{(3)}[2])$ and $(B^{(2)}[1],B^{(3)}[1])$.  Recall
   that each of these pairs has minimum degree at least
   $(1/2-\epsilon^{1/2})t_1$.  Utilizing a lemma in~\cite{Zhao} --
   stated as Lemma~\ref{lem:Zhao} in section~\ref{sec:appfigTWO} --
   we are able to find such a factor unless at least one of those
   pairs is $\alpha(\epsilon^{1/2})$-approximately
   $\Theta_{2\times 2}(t_1/2)$. (Minimality gives the other
   sparse pair.)

   Lemma~\ref{lem:randpairs} says that if random selections give a
   graph that is approximately $\Theta_{2\times 2}$, then the
   original graph was, too.  So, along with Lemma~\ref{lem:Zhao}, it
   establishes that if, after moving our vertices, we are unable
   to complete our $K_{h,h}$-cover in $(B_i(k),B_j(k))$ with
   nontrivial probability, then the pair $(B_i,B_j)$ is
   $\epsilon'$-approximately $\Theta_{2\times 2}(t_1)$, where
   $\epsilon'=\beta(\alpha(\epsilon^{1/2}))$.

   Since none of the pairs is $\epsilon'$-approximately
   $\Theta_{2\times 2}(t_1)$, we can find the required factor
   of $\left(B^{(1)},B^{(2)},B^{(3)}\right)$ by copies of $K_{h,h}$.

\subsubsubsection{Case 2: Exactly one pair is $\Theta_{2\times
2}(t_1)$}

   Here, we will assume that $B^{(1)}=\widetilde{B}^{(1)}+\widehat{B}^{(1)}$ and
   $B^{(2)}=\widetilde{B}^{(2)}+\widehat{B}^{(2)}$, where
   $|\widetilde{B}^{(1)}|=|\widehat{B}^{(2)}|=t_1$ and
   $d(\widetilde{B}^{(1)},\widehat{B}^{(2)}),d(\widehat{B}^{(1)},\widetilde{B}^{(2)})\leq\epsilon'$.
   A random partition of $B^{(1)}$ into pieces, with probability approaching 1
   as $t_1$ approaches infinity, will partition $\widetilde{B}^{(1)}$ into two
   approximately equal pieces.  In particular, let the {\bf
   typical vertices} in $\widetilde{B}^{(1)}$ be those that are nonadjacent to at
   most $(\epsilon')^{1/2}t_1$ in $\widehat{B}^{(2)}$.  There are at most
   $(\epsilon')^{1/2}t_1$ such vertices.  A similar conclusion can
   be drawn from $\widetilde{B}^{(2)}$, $\widehat{B}^{(1)}$ and $\widehat{B}^{(2)}$.

   In this case, we randomly partition $B^{(1)}$, $B^{(2)}$ and $B^{(3)}$ into
   the sets $B^{(i)}[k]$ as proscribed.  Exchange the vertices as we
   have done above and complete both the reserved and exchanged
   vertices to form copies of $K_{h,h}$.  This encompasses at most
   $8h\epsilon t_1$ vertices.
   Exchange vertices in $B^{(1)}[3]$ with vertices in $B^{(1)}[2]$ and
   vertices in $B^{(2)}[3]$ with vertices in $B^{(2)}[1]$ so that there
   are exactly $h\lfloor t_1/(2h)\rfloor$ typical vertices of
   $\widetilde{B}^{(1)}$ in $B^{(1)}[3]$ and $h\lfloor t_1/(2h)\rfloor$ typical
   vertices of $\widehat{B}^{(2)}$ in $B^{(2)}[3]$.  Let the rest of the vertices,
   not matched into a $K_{h,h}$, in $B^{(1)}[3]$ be typical vertices
   in $\widehat{B}^{(1)}$ and the rest of the
   vertices in $B^{(2)}[3]$ be typical in $\widetilde{B}^{(2)}$.  Using
   Proposition~\ref{prop:sufffactor}(\ref{it:bi:sufffactor}) on
   each pair of sets of typical vertices in $(B^{(1)}[3],B^{(2)}[3])$
   will easily have a $K_{h,h}$-factor.  With $\epsilon'$ small enough, we can
   guarantee that at most $(\epsilon')^{1/3}t_1$ vertices in
   $(B^{(1)}[2],B^{(3)}[2])$ and $(B^{(2)}[1],B_3[1])$ were moved.  Applying
   Lemmas~\ref{lem:Zhao} and~\ref{lem:randpairs}, and the fact that no
   pair other than $(B^{(1)},B^{(2)})$ can be $\epsilon'$-approximately
   $\Theta_{2\times 2}(t_1)$, we conclude that the pairs
   $(B^{(1)}[2],B^{(3)}[2])$ and $(B^{(2)}[1],B_3[1])$
   can be completed to $K_{h,h}$-factors.

\subsubsubsection{Case 3: Exactly two pairs are $\Theta_{2\times
2}(t_1)$, which do not coincide}

   Let the pairs in question be $(B^{(1)},B^{(2)})$ and $(B^{(2)},B^{(3)})$.
   Let the dense pairs in the subgraph induced by $(B^{(1)},B^{(2)})$ be
   $(\widetilde{B}^{(1)},\widetilde{B}^{(2)})$ and
   $(\widehat{B}^{(1)},\widehat{B}^{(2)})$.  Let the dense pairs in
   $(B^{(2)},B^{(3)})$ be $(\mathring{B}^{(2)},\mathring{B}^{(3)})$ and
   $(\ddot{B}^{(2)},\ddot{B}^{(3)})$.  Moreover, since the pairs
   fail to coincide, we can conclude that the intersection of the typical
   vertices of $\widetilde{B}^{(2)}$ with the typical vertices of each of
   $\mathring{B}^{(2)}$ and $\ddot{B}^{(2)}$ is at least
   $(\epsilon')^{1/4}t_1$ and similarly for $\widehat{B}^{(2)}$.

   Once again, we randomly partition the vertices in $B^{(1)}$, $B^{(2)}$
   and $B^{(3)}$ and move vertices so as to ensure that the reserved
   vertices and the vertices exchanged for them are placed
   into vertex-disjoint copies of $K_{h,h}$.  Our concern at this
   point is the vertices in $B^{(2)}$.

   Consider the vertices in $(B^{(1)}[3],B^{(2)}[3])$.
   Approximately half are typical vertices of $\widetilde{B}^{(2)}$
   and approximately half are typical vertices of
   $\widehat{B}^{(2)}$.  Take each non-typical vertex in $B^{(1)}[3]$
   and in $B^{(2)}[3]$, match them with a copy of $K_{h,h}$ in the
   pair $(B^{(1)}[3],B^{(2)}[3])$ and remove them.  Do the same for
   vertices in $B^{(2)}[1]$ that are not typical in
   $\mathring{B}^{(2)}$ or $\ddot{B}^{(2)}$ and in $B^{(3)}[1]$ that
   are not typical in $\mathring{B}^{(3)}$ or $\ddot{B}^{(3)}$.
   Remove those copies of $K_{h,h}$ also.

   Observe that there are at least $\epsilon^{1/4}t_1/4$ vertices
   in each intersection of $\widetilde{B}^{(2)}$ or $\widehat{B}^{(2)}$
   with $\mathring{B}^{(2)}$ or $\ddot{B}^{(2)}$ and with
   $B^{(2)}[3]$ or $B^{(2)}[1]$.

   First, move $a$ vertices from $\widetilde{B}^{(2)}\cap
   \mathring{B}^{(2)}\cap B^{(2)}[3]$ to $\widetilde{B}^{(2)}\cap
   \mathring{B}^{(2)}\cap B^{(2)}[1]$ to make $|\widetilde{B}^{(2)}\cap
   B^{(2)}[3]|$ divisible by $h$.  Second, move $a+b$ vertices from
   $\widehat{B}^{(2)}\cap \mathring{B}^{(2)}\cap B^{(2)}[1]$ to
   $\widehat{B}^{(2)}\cap \mathring{B}^{(2)}\cap B^{(2)}[3]$ to
   make $|\mathring{B}^{(2)}\cap B^{(2)}[1]|$ divisible by $h$.  Third,
   move $a+b+c$ vertices from $\widehat{B}^{(2)}\cap \ddot{B}^{(2)}\cap
   B^{(2)}[3]$ to $\widehat{B}^{(2)}\cap \ddot{B}^{(2)}\cap B^{(2)}[1]$.
   This will make both $|\widehat{B}^{(2)}\cap B^{(2)}[3]|$ and
   $|\ddot{B}^{(2)}\cap B^{(2)}[1]|$ divisible by $h$.

   Here $a$, $b$ and $c$ are the remainders of $|\widetilde{B}^{(2)}\cap
   B^{(2)}[3]|$, $|\mathring{B}^{(2)}\cap B^{(2)}[1]|$ and
   $|\widehat{B}^{(2)}\cap B^{(2)}[3]|$, respectively, when each is divided by
   $h$.  Observe that both $|\widetilde{B}^{(2)}\cap
   B^{(2)}[3]|+|\widehat{B}^{(2)}\cap B^{(2)}[3]|$ and
   $|\mathring{B}^{(2)}\cap B^{(2)}[1]|+|\ddot{B}^{(2)}\cap
   B^{(2)}[1]|$ are divisible by $h$.

   Finally, we exchange vertices in $\widetilde{B}^{(1)}\cap B^{(1)}[3]$ with
   those in $\widetilde{B}^{(1)}\cap B^{(1)}[2]$ so that $|\widetilde{B}^{(1)}\cap
   B^{(1)}[3]|=|\widetilde{B}^{(2)}\cap B^{(2)}[3]|$ and similarly for
   $\widehat{B}^{(2)}$.  Also, exchange vertices
   in $\mathring{B}^{(3)}\cap B^{(3)}[1]$ with those in
   $\mathring{B}^{(3)}\cap B^{(3)}[2]$ so that
   $|\mathring{B}^{(3)}\cap B^{(3)}[1]|=|\mathring{B}^{(2)}\cap
   B^{(2)}[1]|$ and similarly for $\ddot{B}^{(2)}$.

   Then, in $(\widetilde{B}^{(1)}\cap B^{(1)}[3],\widetilde{B}^{(2)}\cap
   B^{(2)}[3])$, first greedily place each moved vertex into copies of
   $K_{h,h}$ and then finish the factor via
   Proposition~\ref{prop:sufffactor}(\ref{it:bi:sufffactor}).
   Do the same for $\left(\widehat{B}^{(1)}\cap B^{(1)}[3],\widehat{B}^{(2)}\cap
   B^{(2)}[3]\right)$,
   $\left(\mathring{B}^{(2)}\cap B^{(2)}[1],\mathring{B}^{(3)}\cap
   B^{(3)}[1]\right)$ and $\left(\mathring{B}^{(2)}\cap
   B^{(2)}[1],\mathring{B}^{(3)}\cap B^{(3)}[1]\right)$.

   Finally, we can complete the factor of $(B^{(1)}[2],B^{(3)}[2])$
   because if it is not possible, Lemmas~\ref{lem:Zhao}
   and~\ref{lem:randpairs} would require $(B^{(1)},B^{(3)})$ to be
   approximately $\Theta_{2\times 2}(t_1)$, excluded by this case.

\subsubsubsection{Case 4: Three pairs are $\Theta_{2\times 2}(t_1)$,
none of which coincide}

   Let the dense pairs in $(B^{(1)},B^{(2)})$ be
   $(\widetilde{B}^{(1)},\widetilde{B}^{(2)})$ and
   $(\widehat{B}^{(1)},\widehat{B}^{(2)})$.
   Let the dense pairs in $(B^{(2)},B^{(3)})$ be
   $(\mathring{B}^{(2)},\mathring{B}^{(3)})$ and
   $(\ddot{B}^{(2)},\ddot{B}^{(3)})$.  Let the dense
   pairs in $(B^{(1)},B^{(3)})$ be
   $(B^{(1)}_{\sharp},B^{(3)}_{\sharp})$ and
   $(B^{(1)}_{\flat},B^{(3)}_{\flat})$.  Moreover,
   since the pairs fail to coincide, we can conclude that the
   intersection of the typical vertices of one set of sparse pairs
   with the typical vertices of another is at least
   $(\epsilon')^{1/4}t^{(1)}$.

   Partition $B^{(1)}$, $B^{(2)}$ and $B^{(3)}$ into
   appropriately-sized
   sets as before, uniformly at random.  The degree conditions hold
   with high probability as before.  Take non-typical vertices and
   complete them greedily to place them in vertex-disjoint copies
   of $K_{h,h}$ within each of the pairs $(B^{(1)}[3],B^{(2)}[3])$,
   $(B^{(2)}[1],B^{(3)}[1])$ and $(B^{(1)}[2],B^{(3)}[2])$.  Remove
   these
   copies of $K_{h,h}$ from the graph.

   Let $M$ be the largest multiple of $h$ less than or equal to
   the size of the intersection of what remains of any sparse set
   (\ie, $\widetilde{B}^{(i)},\widehat{B}^{(i)},
   \mathring{B}^{(i)},\ddot{B}^{(i)},
   B^{(i)}_{\sharp},B^{(i)}_{\flat}$) with a set of the form
   $B^{(i)}[k]$.

   We can move vertices as in Case 3 by letting
   $a=|\widetilde{B}^{(2)}\cap B^{(2)}[3]|-M$,
   $b=|\mathring{B}^{(2)}\cap B^{(2)}[1]|-M$ and
   $c=|\widehat{B}^{(2)}\cap B^{(2)}[3]|+M-t_3$, which is also equal
   to $t_1-M-a-b-|\ddot{B}^{(2)}\cap B^{(2)}[1]|$.  We can
   perform similar operations to guarantee that, among the vertices
   that remain in the graph, that
   \begin{eqnarray*}
   M & = & \left|\widetilde{B}^{(1)}\cap
   B^{(1)}[3]\right|=\left|\widetilde{B}^{(2)}\cap B^{(2)}[3]\right|=\left|\mathring{B}^{(2)}\cap
   B^{(2)}[1]\right|=\left|\mathring{B}^{(3)}\cap
   B^{(3)}[1]\right| \\
   & = & \left|B^{(1)}_{\sharp}\cap
   B^{(1)}[2]\right|=\left|B^{(3)}_{\sharp}\cap B^{(3)}[2]\right|
   \end{eqnarray*}
   The fact that the pairs do not coincide ensures that there are
   enough vertices to make these moves.

   Place the moved vertices into vertex-disjoint copies of
   $K_{h,h}$ and finish the factor via
   Proposition~\ref{prop:sufffactor}(\ref{it:bi:sufffactor}).
\end{proofcite}

\begin{proofcite}{Proposition~\ref{prop:sufffactor}}
   \begin{enumerate}
      \item This is found by arbitrarily placing vertices
      from the same part into clusters of size $h$.  Construct
      an auxiliary graph $G'$ on the clusters where two are
      adjacent if and only if they form a $K_{h,h}$ in $G$.
      Each cluster in $G'$ is adjacent to at least half of the
      $M/h$ clusters in the other part. Using K\"onig-Hall, we
      find a matching in $G'$, producing a $K_{h,h}$-factor.
      \item The idea is the same as above -- place vertices into
      clusters of size $h$ -- and use the tripartite version of
      Proposition~1.3 in~\cite{MSz} as a generalization of
      K\"onig-Hall.
   \end{enumerate}
\end{proofcite}

\section{Lower bounds}
\label{sec:lb}

We give a number of constructions which establish the lower bounds. The constructions in~\cite{Zhao} of sparse regular bipartite graphs with no $C_4$'s lead naturally to the following
important proposition, which we state without proof.
\begin{proposition}
   For each integer $d\geq 0$, there exists an $n_0$ such that, if
   $n\geq n_0$, there exists a balanced tripartite graph, $Q(n,d)$,
   on $3n$ vertices such that each of the ${3\choose 2}$ natural
   bipartite subgraphs are $d$-regular with no $C_4$ and $Q(n,d)$
   has no $K_3$.
\end{proposition}

\subsection{Tight lower bound for $(6h)\mid N$}
Recall that if $G\in{\cal G}_3(N)$, $N\geq N_0$ has minimum degree at least $h\left\lceil\frac{2N}{3h}\right\rceil+(h-1)$ and is not in~\vexc, then $G$ has a $K_{h,h,h}$-factor.  Proposition~\ref{prop:LB0} shows that our results are best possible in the case where $N$ is a multiple of $6h$ or even $N$ is a multiple of $3h$ but the graph is not in~\vexc.
\begin{proposition}
   Fix a natural number $h\geq 2$ and $N=3qh$.  If $q$ is large enough, there exists a $G_0\in{\cal G}_3(N)$ such that $\bar{\delta}(G_0)=h\left\lceil\frac{2N}{3h}\right\rceil+h-2=2qh+(h-2)$
   and $G_0$ has no $K_{h,h,h}$-factor.
   \label{prop:LB0}
\end{proposition}

\begin{proof}
   We will construct 9 sets $A^{(i)}_{j}$ with $i,j\in\{1,2,3\}$.
   The union $A^{(i)}_{1}+A^{(i)}_2+A^{(i)}_3$ defines the
   \textbf{$i^{\rm th}$ vertex-class}.
   Call the triple $(A^{(1)}_{j},A^{(2)}_{j},A^{(3)}_{j})$ the
   \textbf{$j^{\rm th}$ column}.

   Construct $G_0$ as follows: For $i=1,2,3$, let $|A^{(i)}_{1}|=qh-1$, $|A^{(i)}_{2}|=qh$
   and $|A^{(i)}_{3}|=qh+1$.  Let the graph in column 1 be
   $Q(qh-1,h-3)$, the graph in column 2 be
   $Q(qh,h-2)$ and the graph in column 3 be
   $Q(qh+1,h-1)$.  If two vertices are in different columns and different vertex-classes, then they are adjacent.  It is easy to verify that $\bar{\delta}(G_0)=2qh+(h-2)$.  Suppose, by way of contradiction, that $G_0$ has a $K_{h,h,h}$-factor.

   Since there are no triangles and no $C_4$'s in any column, the intersection of a copy of $K_{h,h,h}$ with a column is either a star, with all leaves in the same vertex-class, or a set of vertices in the same vertex-class.  So, each copy of $K_{h,h,h}$ have at most $h$ vertices in column 3.  A $K_{h,h,h}$-factor has exactly $3q$ copies of $K_{h,h,h}$ and so the factor has at most $3qh$ vertices in column 3.  But there are $3qh+3$ vertices in column 3, a contradiction.
\end{proof}

\subsection{General lower bound for $h\mid N$}
Proposition~\ref{prop:LB12} gives a more general lower bound for cases when $N/h$ is not divisible by $3$, although it leaves a gap of 1 from the upper bound.
\begin{proposition}
   Fix a natural number $h\geq 2$ and $N=(3q+r)h$ for $r\in\{0,1,2\}$.  If $q$ is large enough, there exists a $G_1\in{\cal G}_3(N)$ such that $\bar{\delta}(G_1)=h\left\lceil\frac{2N}{3h}\right\rceil+h-3=2qh+rh+(h-3)$
   and $G_1$ has no $K_{h,h,h}$-factor.
   \label{prop:LB12}
\end{proposition}

\begin{proof}
   Define $G_1$ as follows: For $i=1,2,3$, let $|A^{(i)}_{1}|=qh+rh-1$, $|A^{(i)}_{2}|=qh$
   and $|A^{(i)}_{3}|=qh+1$.  Let the graph in column 1 be
   $Q(qh+rh-1,rh+h-4)$ if $rh+h-4\geq 0$ and empty otherwise, the graph in column 2 be
   $Q(qh,h-3)$ and the graph in column 3 be
   $Q(qh+1,h-2)$.  If two vertices are in different columns and different vertex-classes, then they are adjacent.  It is easy to verify that $\bar{\delta}(G_1)=2qh+rh+(h-3)$.  Suppose, by way of contradiction, that $G_1$ has a $K_{h,h,h}$-factor.

   Since there are no triangles and no $C_4$'s in any column, the intersection of a copy of $K_{h,h,h}$ with a column is either a star, with all leaves in the same vertex-class, or a set of vertices in the same vertex-class.  So each copy of $K_{h,h,h}$ has at most $h+1$ vertices in column 1, $h$ vertices in column 2 and at most $h$ vertices in column 3.

   There are three cases for a copy of $K_{h,h,h}$.  Case 1 has $h$ vertices in each column.  Case 2 has $h+1$ vertices in column 1, $h-1$ vertices in column 2 and $h$ vertices in column 3.  Case 3 has $h+1$ vertices in column 1, $h$ vertices in column 2 and $h-1$ vertices in column 3.

   Since a $K_{h,h,h}$ having $h$ vertices in column 3 implies the vertices have the same vertex-class, cases 1 and 2 imply that all vertices in column 3 are in the same vertex-class.  Consider case 3.  Having $h$ vertices in column 2 means that all are in the same vertex-class.  Since $h+1$ vertices in column 1 means that they form a star, the remaining $h-1$ vertices in column 3 must be in the same vertex-class (the same vertex-class as the center of the star).  Hence, every copy of $K_{h,h,h}$ has all of its column 3 vertices in the same vertex-class.  Therefore, the number of copies of $K_{h,h,h}$ in a factor is at least $3\left\lceil\frac{qh+1}{h}\right\rceil=3q+3$, a contradiction because the factor has exactly $3q+r\leq 3q+2$ copies of $K_{h,h,h}$.
\end{proof}

Note that in the previous proof, column 1 could be $Q(qh+rh-1,rh+h-3)$ and column 2 could be $Q(qh,h-2)$ and the argument does not change.  Unfortunately, this proof does require that column 3 have degree at most $h-2$ between parts.

\subsection{Lower bound for~\vexc}
Proposition~\ref{prop:LBvexc} gives a graph in~\vexc~that has minimum degree $2N/3+h-2$, which is greater than that of Proposition~\ref{prop:LB12} but is still far from the upper bound of $2N/3+2h-2$.
\begin{proposition}
   Fix a natural number $h\geq 2$ and $N=(6q+3)h$.  If $q$ is large enough, there exists a $G_2\in{\cal G}_3(N)$ in~\vexc~such that $\bar{\delta}(G_2)=h\left\lceil\frac{2N}{3h}\right\rceil+h-2=(4q+2)h+h-2$
   and $G_2$ has no $K_{h,h,h}$-factor.
   \label{prop:LBvexc}
\end{proposition}

\begin{proof}
   Construct $G_2$ as follows: For $i=1,2,3$ and $j=1,2,3$, let $|A^{(i)}_{j}|=2qh+h$.  Let $\left(A^{(1)}_{1},A^{(2)}_{1},A^{(3)}_{1}\right)$ be $Q(2qh+h,h-2)$.  For $i=1,2,3$, let each vertex in $A^{(i)}_{1}$ be adjacent to any vertex in $A^{(i')}_{j'}$ whenever $i'\neq i$ and $j'\neq 1$.  For $j=2,3$, let $\left(A^{(1)}_{j},A^{(2)}_{j},A^{(3)}_{j}\right)$ be a complete tripartite graph and for $i'\neq i$, let $(A^{(i)}_{2},A^{(i')}_{3})$ be a $(h-2)$-regular graph with no $C_4$.  It is easy to verify that $\bar{\delta}(G_2)=2qh+rh+(h-2)$. Suppose, by way of contradiction, that $G_2$ has a $K_{h,h,h}$-factor.

   Since there are no triangles and no $C_4$'s in column 1, the intersection of a copy of $K_{h,h,h}$ with column 1 is either a star, with all leaves in the same vertex-class, or a set of vertices in the same vertex-class.  So, its intersection is at most $h$ vertices.  Since the factor has $6q+3$ members and column 1 has $6qh+3h$ total vertices, each member of the factor has exactly $h$ vertices in column 1.  As a result, those vertices are in the same vertex-class.

   So, the intersection of any member of the $K_{h,h,h}$-factor with columns 2 and 3 is a $K_{h,h}$ with $h$ vertices in each of two vertex-classes.  Suppose this $K_{h,h}$ has vertices in different columns.  Suppose further that it has one vertex in $A^{(1)}_{2}$, then there are at most $h-2$ vertices in $A^{(2)}_{3}$.  So, there must be at least 2 vertices in $A^{(2)}_{2}$ and since there are no $C_4$'s, at most 1 vertex in $A^{(1)}_{3}$ if $h\geq 3$.  If $h=2$, then there can be no vertices in $A^{(1)}_{3}$.  Regardless, this is a contradiction to the assumption that a $K_{h,h}$ has vertices in both column 2 and column 3.

   So each member of the $K_{h,h,h}$-factor either has $2h$ vertices in column 2 or $2h$ vertices in column 3.  However, there are at most $\left\lfloor\frac{3(2qh+h)}{2h}\right\rfloor=3q+1$ members of the factor with $2h$ vertices in column 2 and at most $3q+1$ members of the factor with $2h$ vertices in column 3.  In either case, there are less than $h$ vertices in $A_2^{(1)}\cup A_3^{(1)}$, a contradiction.
\end{proof} 

\section{Thanks}
The authors would like to thank the Department of Mathematics, Statistics, and Computer Science at the University of Illinois at Chicago for their supporting the first author via a visitor fund for the purposes of working on this project.

\end{document}